# Forwards attractors for non-autonomous Lotka-Volterra cooperative systems: a detailed geometrical description

Juan Garcia-Fuentes[1], Piotr Kalita[1,2], José A. Langa[1], Antonio Suárez[1]

[1]Departamento de Ecuaciones Diferenciales y Análisis Numérico, Universidad de Sevilla, Campus Reina Mercedes, 41012, Sevilla, Spain

[2]Faculty of Mathematics and Computer Science, Jagiellonian University, ul. Lojasiewicza 6, 30-348 Krakow, Poland

E-mail addresses: J.G.-F.: `jgfuentes@us.es` J.A.L.: `langa@us.es` P. K.: `piotr.kalita@ii.uj.edu.pl` A.S.: `suarez@us.es`

February 8, 2024

## Abstract

Non-autonomous differential equations exhibit a highly intricate dynamics, and various concepts have been introduced to describe their qualitative behavior. In general, it is rare to obtain time dependent invariant compact attracting sets when time goes to plus infinity. Moreover, there are only a few papers in the literature that explore the geometric structure of such sets. In this paper we investigate the long time behaviour of cooperative $n$-dimensional non-autonomous Lotka–Volterra systems is population dynamics. We provide sufficient conditions for the existence of a globally stable (forward in time) entire solution in which one species becomes extinct, or where all species except one become extinct. Furthermore, we obtain the precise geometrical structure of the non-autonomous forward attractor in one, two, and three dimensions by establishing heteroclinic connections between the globally stable solution and the semi-stable solutions in cases of species permanence and extinction. We believe that understanding time-dependent forward attractors paves the way for a comprehensive analysis of both transient and long-term behavior in non-autonomous phenomena.

## 1 Introduction

The analysis of non-autonomous dynamical systems has received significant attention in the past three decades. In particular, two main research branches have been developed to understand the long-time behaviour of such systems: on the one hand, pullback dynamics, which focuses on observing the system state at time $t$ after it has evolved for a significant duration from the past. On the other hand, forwards dynamics focuses on the long time behavior of the system as time tends to plus infinity [CLR13, BCL20, CV02, KR11]. In the former case, we get a time dependent invariant family of sets; in the second one, the object to study, the uniform attractor, is the smallest compact (not necessarily invariant) set attracting uniformly all bounded sets for all initial times. But there exists a very interesting scenario which provides the new insight into the asymptotic behavior of non-autonomous systems. Indeed, sometimes we observe a time depending invariant family of compact sets which is forwards attracting, combining the advantageous properties of both approaches, namely invariance from pullback dynamics and forward attraction from uniform attractor theory. However, obtaining these sets is generally a rare occurrence, and a well-established theory for their existence is still lacking. Moreover, one of the most significant and promising properties of attractors is revealed when we know their geometrical structure, i.e., their invariant subsets made of recurrent points and the connections made by gradient trajectories. This structure is maintained when defining the structural stability of these systems [BCL20]. In this paper we study forward attractors for $n$-dimensional Lotka–Volterra systems and provide the detailed characterization of their internal geometrical structure. We believe that both types of results, which we have obtained, are novel and significant contributions to the existing literature.

The Lotka–Volterra system is a classical model of population dynamics representing the behavior of species in ecosystems. This model consists of the following $n$-dimensional system of differential equations, in which the unknown $u_i$ corresponds to the density of individuals of the $i$-th species

$$u_i' = u_i \left( a_i(t) - b_{ii}(t)u_i - \sum_{\substack{j=1 \\ j \neq i}}^{n} b_{ij}(t)u_j \right) \quad \text{for } i \in \{1, ..., n\}. \tag{1}$$

The function $a_i$ determines the intrinsic growth rate of the species $i$, and $b_{ij}$ represents the interaction between the species $i$ and $j$. In particular $b_{ii}$ stands for the intrinsic competition of the same species $i$. If $b_{ij}$ is negative, then the species $j$ contributes positively to the growth of the species $i$, and we speak about a cooperative system. Alternatively, if $b_{ij}$ is positive, $j$ contributes negatively to the growth of $i$ which means that the species compete with the others.

When these coefficients are constant, the system is autonomous. In such case if the matrix $B = (b_{ij})_{i,j=1}^{n}$ is Volterra–Lyapunov stable (see, for instance, [Tak96]), the full characterization of the asymptotic behavior of the solutions of the system is well known. In particular, it is possible to give conditions on the coefficients $(a_i)_{i=1}^{n}$ which guarantee the existence of a globally asymptotically stable equilibrium, in which all the species are present, or, on the contrary, conditions for this stable equilibrium in which one or more species become extinct. These observations allow to determine the full structure of the global attractor associated to (1) [Tak96, Gue17, HS22, AKL22]. However, the situation is entirely different in the non-autonomous case, where coefficients depend on time. The corresponding theory is not fully developed, particularly concerning the detailed description of the attractor geometric structure. In the present paper we study the non-autonomous cooperative Lotka–Volterra system and we give conditions to describe its long time behavior, including the full characterization of its attractor in low dimensional cases.

We briefly recall the most important results obtained previously for this non-autonomous version of the problem. In the works of Gopalsamy [Gop86a, Gop86b] the author treats the case of periodic and almost periodic functions $a_i(t)$ and $b_{ij}(t)$ in competitive framework, i.e. $b_{ij}(t) \geqslant 0$ for all $t \in \mathbb{R}$ and $i \neq j$. Conditions are given to ensure the existence of globally stable solutions with all components strictly positive for all times $t \in \mathbb{R}$. Ahmad and Lazer extended the study of the competitive case, relaxing the requirement for coefficients to be periodic or almost periodic. They were still able to establish the existence of a unique globally stable strictly positive permanent solution, in which all species are present [AL95, AL00], or with one of the species extinct [Ahm93, AL98]. On the other hand, Redheffer in [Red96] found some general conditions on the coefficients to get a globally stable solution with no restriction on the sign for the interspecific parameters, but only for the case of the permanence of the system.

In the first part of this paper we obtain the result, not considered before, providing criteria for the extinction of any given subcommunity of species in the cooperative case. The main difficulty of the proof was to obtain the decay of the density of species despite the possible cooperation from the other species. In Lemma 4.1, we obtain that under the diagonal dominance condition $(H1)$ for the non-autonomous matrix $B = (b_{ij})_{i,j=1}^{n}$ and condition $(B)$, which links the matrix $B$ with the intrinsic growth rates of individual species $a = (a_i)_{i=1}^{n}$, densities for a given subcommunity of species decay jointly to zero. Condition $(B)$ is a variation of the conditions in [Red96] although we extend them by allowing some of the growth rates $a_i(t)$ to be negative. We also prove the existence of a stable solution for two particular cases of extinction in the $n$-dimensional species community: one case, in which all species except one become extinct and another case, where only one species becomes extinct. These results are contained in Sections 4.3 and 4.4, respectively.

The feature of the system (1) is, that by setting some of the variables to zero, one obtains, for each subset of the community $\{1, \ldots, n\}$, a smaller system, modeling the dynamics of the corresponding subcommunity, embedded in the original system. Hence, the above results of Redheffer, Ahmad and Lazer as well as the ones in our Section 4 allow us to construct the asymptotically stable solutions for the full community as well as all its subcommunities. Here a natural question appears, motivated by recent developments in the theory of attractors for non-autonomous dynamical systems [KR11, CLR13, BCL20]: how to describe all solutions which are bounded both in the past and future, and how are they related to each other in terms of attraction when time tends to plus or minus infinity. This question naturally leads to the study of the geometrical structure of non-autonomous attractors, which consist of the globally asymptotically stable solutions, which were the main interest to previous researchers analysing (1), asymptotically stable solutions for the subcommunities, and heteroclinic connections between them. By the theory of exponential dichotomies for the linearized systems (see Section 2), based on the results of [BP15] on dichotomies in upper triangular systems, and the non-autonomous versions of the Hadamard–Perron theorem on the existence of local stable and unstable manifolds [CLR13, KR11] we are able to construct the heteroclinic connections between various stable and semi-stable (i.e. stable for the subsystems of the original systems) solutions. We also prove that all other solutions must be unbounded in the past. Thus, we show in Sections 5, 6, and 7 the complete structure of the non-autonomous attractor for various cases of extinction and persistence, for one, two, and three dimensions, respectively. Such detailed description of a non-autonomous attractor is always a very difficult task, and we have only a few very recent papers on this subject (such as [CRLO23]). However, thanks to the previous results, we construct the full structure of the forward (and pullback) attractor for the non-autonomous problem governed by the system (1) of differential equations.

# 2 Non-autonomous cooperative Lotka–Volterra system and existence of global solutions

## 2.1 Lotka–Volterra system

We use the following convention for the functions defined on the real line: $a^L = \inf_{t \in \mathbb{R}} a(t)$, $a^U = \sup_{t \in \mathbb{R}} a(t)$. Moreover, for two such functions $a, b$ we will say that $a > b$ if $a(t) > b(t)$ for $t \in \mathbb{R}$ and likewise we define the non-sharp inequalities $a \geqslant b$.

We will consider the following $n$-dimensional cooperative Lotka–Volterra system.

$$u_i' = u_i\left(a_i(t) - \sum_{j=1}^{n} b_{ij}(t)u_j\right), \tag{LV-$n$}$$

for $i = 1, ..., n$, where $a_i(t)$ and $b_{ij}(t)$ are continuous real-valued functions. Assume that the system is cooperative, i.e. $b_{ij}(t) \leqslant 0$ for $i \neq j$, and $b_{ii}(t) > 0$ for all $t \in \mathbb{R}$, and for $i = 1, ..., n$. Moreover, we make the following standing assumption of the row diagonal dominance

$$\left(c_i b_{ii} + \sum_{j=1, j \neq i}^{n} c_j b_{ij}\right)^L \geqslant \delta > 0 \quad \text{for all } i = 1, ..., n, \tag{H1}$$

for some positive constants $\{c_i\}_{i=1}^{n}$, and column diagonal dominance

$$\left(\bar{c}_i b_{ii} + \sum_{j=1, j \neq i}^{n} \bar{c}_j b_{ji}\right)^L \geqslant \delta > 0 \quad \text{for all } i = 1, ..., n, \tag{H2}$$

for another positive constants $\{\bar{c}_i\}_{i=1}^{n}$.

## 2.2 Existence of global attracting solution. Case of permanence

We recall the result in [Red96] which provides conditions for the existence of a positive solution $u^* = (u_1^*, ..., u_n^*)$, isolated from zero and infinity, which is forwards globally attractive for all trajectories with positive initial data $u(t_0) > 0$.

We need the following assumption: there exist positive vectors $d > 0$ and $\bar{d} > 0$ such that

$$\bar{d}_i b_{ii} \leqslant a_i \leqslant d_i b_{ii} + \sum_{j \neq i} d_j b_{ij} \quad \text{for all } i = 1, ..., n. \tag{A}$$

The following result comes from Redheffer ([Red96, Theorem 1 (ii), (v), and (vi)]):

**Theorem 2.1.** *Assume $(A)$ and $(H2)$. Then there exists a unique complete positive trajectory $u^*$ separated from zero and infinity. This trajectory satisfies $\bar{d}_i \leqslant u_i^* \leqslant d_i$ for every $t \in \mathbb{R}$. Moreover, for every solution $u$ of* (LV-$n$) *such that $u(t_0) > 0$ the following convergence holds*

$$\lim_{t \to \infty} |u(t) - u^*(t)| \to 0.$$

The following result establishes Lipschitz continuous dependence on initial data for forward bounded solutions.

**Lemma 2.2.** *Let two solutions $u$ and $v$ with the initial data at $t_0$ be bounded for $t \geqslant t_0$. Assume that $|a_i|, |b_{ij}|$ are bounded by a constant independent of time. Then there exists $0 < \kappa$ (depending on the bounds on the coefficients and both solutions) such that*

$$|u(t_0 + t) - v(t_0 + t)| \leqslant e^{\kappa t}|u(t_0) - v(t_0)| \quad for \ \ t \geqslant 0.$$

*Proof.* We skip the proof which follows in a standard way from the Gronwall lemma and estimation of the difference of two solutions. □

If the solutions are separated from zero, a stronger result holds under the column diagonal dominance assumption $(H2)$. The following lemma is proved in [Red96, Lemma 8].

**Lemma 2.3.** *Assume $(H2)$ and let $u, v$ be two solutions of (LV-$n$) such that there exist positive constants $\sigma_1, \sigma_2 > 0$ with $\sigma_1 \leqslant u_i(\tau), v_i(\tau) \leqslant \sigma_2$ for every $\tau \geqslant t_0$ and $i \in \{1, \dots, n\}$. Then*

$$|u(t_0 + t) - v(t_0 + t)| \leqslant \frac{\sigma_2}{\sigma_1}|u(t_0) - v(t_0)|e^{-\delta\sigma_1 t},$$

*for every $t \geqslant 0$.*

## 2.3 Existence of global solution. Case of extinction

We provide a generalized version of the above condition $(A)$ for which there exists the global solution consisting of zeros on a subset of indexes denoted by $J$ and far from zero and infinity on the remaining subset of indexes, denoted by $I$. Namely, we replace $(A)$ with the following hypothesis

$$0 < \bar{d}_i b_{ii} \leqslant a_i \leqslant d_i b_{ii} + \sum_{j \neq i, j \in I} d_j b_{ij} \quad \text{for all } i \in I \subset \{1, \dots, n\}, \tag{$A_I$}$$

for positive vectors $(d_j)_{j \in I}$ and $(\bar{d}_j)_{j \in I}$.

As an application of Theorem 2.1 we have the next result. If the indexes are subdivided into two disjoint sets $I$ and $J$, it is possible to assume that zero is a solution on the set $J$. Observe that $(H2)$ holds for any subset of indexes.

**Theorem 2.4.** *If $(H2)$ and $(A_I)$ hold, then there exists a complete trajectory $u^*$ of the (LV-$n$) such that $u_i^* \equiv 0$ for $i \in J = \{1, \dots, n\} \backslash I$ and $\bar{d}_i \leqslant u_i^* \leqslant d_i$ for $i \in I$. For every solution $u$ of (LV-$n$) such that $u_i(t_0) > 0$ for $i \in I$ and $u_i(t_0) = 0$ for $i \in J$ the following convergence holds*

$$\lim_{t \to \infty} |u(t) - u^*(t)| \to 0.$$

**Remark 1.** *If $J$ is nonempty we call the solution semitrivial, as it is equal to zero on some coordinates. Observe that if $(A)$ holds then $(A_I)$ holds for every set of indexes $I$. In the same way, if $(A_I)$ holds, then the same condition is satisfied on any subset of $I$. Note that while Theorem 2.1 gives the convergence of solutions with nonzero initial data for all coordinates, in above theorem we only have convergence for nonzero initial data on coordinates indexed by $I$. We will provide conditions which guarantee that solutions with all nonzero initial data also converge towards the semitrivial solution with some of the species (indexed by $J$) extinct.*

# 3 Exponential dichotomies and linearization

## 3.1 Exponential dichotomies for upper triangular systems

We will consider the following linear nonautonomous system of ODE's

$$x'(t) = D(t)x(t), \tag{2}$$

where $D : \mathbb{R} \to \mathbb{R}^{n \times n}$ is continuous and bounded, i.e. $\sup_{t \in \mathbb{R}} \|D(t)\| \leqslant M$, where $M$ is a positive constant. The solutions of this system of ODE's, with the initial condition $x(t_0) = x_0$ are given by the multiplication of the initial data by the following fundamental matrix $x(t) = M_D(t, t_0)x_0$. This fundamental matrix is invertible and if $t < t_0$, then $M_D(t, t_0) = M_D(t_0, t)^{-1}$. When $t_0 = 0$, we use the notation $M_D(t, 0) = M_D(t)$. Then $M_D(t, s) = M_D(t)M_D(s)^{-1}$. Moreover the fundamental matrix satisfies the following variational ODE $M_D'(t) = D(t)M_D(t)$ for every $t \in \mathbb{R}$ with the initial data $M_D(0) = I$.

**Definition 3.1.** *We say that the system* (2) *has an exponential dichotomy on an interval $I$ with projection $P : I \to \mathbb{R}^{n \times n}$, constant $k \geqslant 1$ and exponents $\alpha, \beta > 0$ if the fundamental matrix satisfies the invariance property*

$$P(t)M_D(t, s) = M_D(t, s)P(s) \quad \text{for all} \ \ t, s \in I, \tag{3}$$

*and we have the inequalities*

$$\|M_D(t, s)P(s)\| \leqslant ke^{-\alpha(t-s)} \quad \text{for all} \ \ s \leqslant t \in I, \tag{4}$$

$$\|M_D(t, s)(I - P(s))\| \leqslant ke^{\beta(t-s)} \quad \text{for all} \ \ t \leqslant s \in I. \tag{5}$$

Typically $I = \mathbb{R}$, $I = [0, \infty)$ or $I = (-\infty, 0]$. In each of these cases if we substitute $s = 0$ in (3) we obtain $P(t)M_D(t) = M_D(t)P(0)$, and hence $P(t) = M_D(t)P(0)M_D(t)^{-1}$, so we can recover $P(t)$ from $P(0)$ and the fundamental matrix. We denote $P(0) = P$. Substituting $s = t$ in (4) we obtain $\|P(t)\| \leqslant k$. This means that if system (2) has an exponential dichotomy, then $\|P(t)\|$ is bounded uniformly with respect to $t \in I$, and hence the moduli of all entries of the matrix $P(t)$ must be also bounded. If $I = \mathbb{R}$ then the projection $P$ is given uniquely, cf. [Cop78, page 19]. If $I = \mathbb{R}^-$, then $P$ in the definition of the dichotomy can be replaced by any other projection with the same kernel, cf. [BF20, Proposition 2], and if $I = \mathbb{R}^+$, then $P$ can be replaced by any other projection with the same range. If the equation has the dichotomy on $\mathbb{R}^+$ and $\mathbb{R}^-$ with the same projection $P$ and exponents then it also has the dichotomy on $\mathbb{R}$ with the same projection and exponents, cf. [BF20, Proposition 1] or [Cop78, page 19].

We recall the results of [BP15] on upper triangular systems. We will consider the problem governed by the system

$$\begin{pmatrix} x(t) \\ y(t) \end{pmatrix}' = \begin{pmatrix} A(t) & C(t) \\ 0 & B(t) \end{pmatrix} \begin{pmatrix} x(t) \\ y(t) \end{pmatrix}, \tag{6}$$

assuming that $A(t), B(t), C(t)$ are bounded and continuous. We discuss the relation between the fact that the smaller systems governed by the matrices given by the diagonal blocks $x'(t) = A(t)x(t)$ and $y'(t) = B(t)y(t)$ have the dichotomies, with the fact that (6) has the exponential dichotomy. We cite the following result, cf [BP15, Corollary 1]

**Corollary 3.2.** *Assume that the linear systems $x'(t) = A(t)x(t)$ and $y'(t) = B(t)y(t)$ have exponential dichotomies on $\mathbb{R}$, where $A(t) \in \mathbb{R}^{d \times d}$ and $B(t) \in \mathbb{R}^{(n-d) \times (n-d)}$, and $C(t)$ is piecewise continuous and bounded $d \times (n-d)$ matrix. Then (6) has the exponential dichotomy on $\mathbb{R}$.*

The projection $P$ for (6) can be constructed in the following way. First we construct the projections for the dichotomies on $\mathbb{R}^-$ and $\mathbb{R}^+$ in the following way, cf. [BP15, Theorem 1]

$$P^+ = \begin{pmatrix} P^A & L^+ P^B \\ 0 & P^B \end{pmatrix} \text{ on } \mathbb{R}^+ \qquad P^- = \begin{pmatrix} P^A & L^-(I_{n-d} - P^B) \\ 0 & P^B \end{pmatrix} \text{ on } \mathbb{R}^-, \tag{7}$$

where $P^A$ and $P^B$ are the projections for systems with $A(t)$ and $B(t)$ respectively, and $L^+$, $L^-$ are the linking operators given by

$$L^+ = -\int_0^\infty M_A(s)^{-1}(I_d - P^A(s))C(s)M_B(s)\,ds,$$

$$L^- = \int_{-\infty}^0 M_A(s)P^A(s)C(s)M_B(s)^{-1}\,ds.$$

Next, we construct the unique projection $P$, such that its kernel coincides with the kernel of $P^-$ and its range with the range of $P^+$. This projection gives the dichotomy on $\mathbb{R}$.

## 3.2   Linearization of nonautonomous Lotka–Volterra system

Assume that $(u_1^*(t), u_2^*(t), \ldots, u_k^*(t), 0, \ldots, 0)$ is a solution of (LV-$n$), such that $u_i^*$ are separated from zero and infinity for $i \in \{1, \ldots, k\}$. We consider the linearized system around this solution. Denoting $u_j^* = 0$ for $j \in \{k+1, \ldots, n\}$ and $w_i = u_i - u_i^*$ we obtain the system of differential equations

$$w_i'(t) = a_i(t)w_i(t) - u_i^*(t)\sum_{j=1}^n b_{ij}(t)w_j(t) - w_i(t)\sum_{j=1}^k b_{ij}(t)u_j^*(t) - \sum_{j=1}^n b_{ij}(t)w_j(t)w_i(t). \tag{8}$$

The above system can be written as $w'(t) = \mathcal{M}(t)w(t) + R(w(t), t)$, where $R$ is a remainder, quadratic with respect to $w$. Dropping this quadratic term we get the linearized system $v'(t) = \mathcal{M}(t)v(t)$, where the unknown is now denoted by $v$. The system has the following block diagonal form

$$v'(t) = \begin{pmatrix} A(t) & C(t) \\ 0 & B(t) \end{pmatrix} v(t), \tag{9}$$

where $A(t)$ is $k \times k$ matrix with entries given by $a_{ii}(t) = a_i(t) - \sum_{j=1}^k b_{ij}(t)u_j^*(t) - b_{ii}(t)u_i^*(t)$ and $a_{ij}(t) = -b_{ij}(t)u_i^*(t)$ for $j \neq i$. The matrix $B(t)$ is an $(n-k) \times (n-k)$ diagonal matrix with $c_{ii}(t) = a_i(t) - \sum_{j=1}^k b_{ij}(t)u_j^*(t)$. Finally, $C(t)$ in an $k \times (n-k)$ matrix with entries given by $c_{ij}(t) = -u_i^*(t)b_{ij}(t)$.

The following Lemma, which has been proved in [AL98, Lemma 3.6], states that the system with the matrix $A(t)$ always has the exponential dichotomy with $P = I$.

**Lemma 3.3.** *Assume the column diagonal dominance condition (H2). Let $u^*(t)$ be a positive solution of (LV-$n$) separated from zero and infinity, and let $A(t) = (a_{ij}(t))$ be an $n \times n$ matrix* *of the system linearized around $u^*(t)$, that is $a_{ii}(t) = a_i(t) - \sum_{j=1}^n b_{ij}(t)u_j^*(t) - b_{ii}(t)u_i^*(t)$ and $a_{ij}(t) = -b_{ij}(t)u_i^*(t)$ for $j \neq i$. Let be $M_A(t)$ will be a fundamental matrix of the system $v'(t) = A(t)v(t)$, that is*

$$M_A'(t) = A(t)M_A(t), \quad M_A(0) = I.$$

*There exist constants $K > 0$ and $\gamma > 0$ such that*

$$\|M_A(t, s)\| = \|M_A(t)M_A^{-1}(s)\| \leqslant Ke^{-\gamma(t-s)}$$

*for $-\infty < s \leqslant t < \infty$.*

Of course the above lemma remains valid if the dimension $n$ is replaced with smaller numer $k$, then the linear stability holds in the $k$ dimensional subspace.

# 4 Asymptotic stability for solutions with extinction

## 4.1 Extinction conditions for a given subset of species

In this section we present the assumptions on $a_i$ and $b_{ij}$ which enable the possibility of splitting the set $\{1, \ldots, n\}$ into the sum of two disjoint sets of indexes $I$ and $J$. For any initial data $u(t_0) > 0$ the species indexed by $J$ go to extinction, that is, $u_i(t) \to 0$ when $t \to \infty$ for $i \in J$ and $u_i(t)$ is separated from zero forward in time when $i \in I$. We impose the following condition.

$$\begin{cases} b_{ii}\bar{d}_i + \varepsilon \leqslant a_i \leqslant b_{ii}d_i + \sum_{j=1, j\neq i}^n b_{ij}(d_j + \theta c_j) - \varepsilon \text{ for } i \in I \\ a_i \leqslant \sum_{j=1, j\neq i}^n b_{ij}(d_j + \theta c_j) - \varepsilon \text{ for } i \in J \end{cases} \quad (B)$$

for some positive vectors $d, \bar{d} > 0$, and positive numbers $\varepsilon, \theta > 0$. The constants $c_j$ appear in $(H1)$. Note that if $i \in J$, then $a_i^U \leqslant -\varepsilon$, and if $i \in I$ then $a_i^L \geqslant \varepsilon$. Hence $a_i$ is negative and separated from zero for $i \in J$ and positive and separated from zero for $i \in I$.

The proof of the following result is the generalization to $n$ dimensions of the result of Ahmad [Ahm93, Lemma 2], where the two dimensional problem is considered and both sets $I$ and $J$ consist of a single index.

**Lemma 4.1.** *If $(B)$ and $(H1)$ hold then for every solution $u$ of (LV-$n$) such that $u(t_0) > 0$ there exists $t^* > t_0$ such that for every $t \geqslant t^*$*

- $\bar{d}_i < u_i(t) < d_i$ for $i \in I$,

- $u_i(t) < d_i$ for $i \in J$.

*Proof.* **Step 1.** We first prove that for every $i \in I$ there exists $t_i = t(i)$ such that $u_i(t) > \bar{d}_i$ for every $t \geqslant t_i$. First suppose that there exists $\underline{t} \geqslant t_0$ such that $u_i(t) > \bar{d}_i$. We prove that $u_i(s) > \bar{d}_i$ for every $s \geqslant t$. Indeed assume that $u_i(s) \leqslant \bar{d}_i$ for some $s > t$. Define $r = \min\{r \geqslant t : u_i(r) = \bar{d}_i\}$. On one hand it must be $\dot{u}_i(r) \leqslant 0$. On the other hand

$$u_i'(r) = \bar{d}_i \left( a_i(r) - b_{ii}(r)\bar{d}_i - \sum_{j=1, j\neq i}^n b_{ij}(r)u_j(r) \right) \geqslant \bar{d}_i \varepsilon,$$

a contradiction. Now assume that there exists $t > t_0$ such that $u_i(s) \leqslant \bar{d}_i$ for every $s \geqslant t$. Then we would have

$$\frac{u_i'(s)}{u_i(s)} = a_i(s) - b_{ii}(s)u_i(s) - \sum_{j=1, j \neq i}^{n} b_{ij}(s)u_j(s) > a_i(s) - b_{ii}(s)\bar{d}_i \geqslant \varepsilon > 0$$

for every $s \geqslant t$. This implies that $u_i(s) \to \infty$ as $s \to \infty$, so we have a contradiction.

**Step 2.** Define the sets

$$A_k = [0, d_1 + k\theta c_1] \times [0, d_2 + k\theta c_2] \times ... \times [0, d_n + k\theta c_n] \quad \text{for} \quad k \geqslant 0. \tag{10}$$

We first prove that if $u_0 \in A_k$, then $u(t) \in A_k$ for all $t \geqslant t_0$, i.e. the sets $A_k$ are positively invariant. Assume that this is not the case, i.e. $u(t_0) \in A_k$ and for some $t \geqslant t_0$ $u(t) \notin A_k$. Then there must exist $s > t_0$ and an index $i$ such that $u_i(s) = d_i + k\theta c_i$ and $\dot{u}_i(s) \geqslant 0$, and $u_j(s) \leqslant d_j + k\theta c_j$ for $j \neq i$. We calculate

$$a_i(s) \geqslant b_{ii}(s)u_i(s) + \sum_{j=1, j \neq i}^{n} b_{ij}(s)u_j(s) = b_{ii}(s)c_i k\theta + b_{ii}(s)d_i + \sum_{j=1, j \neq i}^{n} b_{ij}(s)u_j(s)$$

$$\geqslant b_{ii}(s)c_i k\theta + b_{ii}(s)d_i + \sum_{j=1, j \neq i}^{n} b_{ij}(s)d_j + \sum_{j=1, j \neq i}^{n} b_{ij}(s)c_j k\theta$$

$$= b_{ii}(s)d_i + \sum_{j=1, j \neq i}^{n} b_{ij}(s)d_j + \theta k \left( \sum_{j=1}^{n} c_j b_{ij}(s) \right)$$

$$> b_{ii}(s)d_i + \sum_{j=1, j \neq i}^{n} b_{ij}(s)d_j > b_{ii}(s)d_i + \sum_{j=1, j \neq i}^{n} b_{ij}(s)(d_j + \theta c_j).$$

Now if $i \in I$, we directly get the contradiction with $(B)$, and if $i \in J$, then the above computation implies that

$$a_i(s) > \sum_{j=1, j \neq i}^{n} b_{ij}(s)(d_j + \theta c_j),$$

and we arrive at the contradiction with $(B)$, too.

**Step 3.** Now let $k \geqslant 0$. We prove that if $u(t) \in A_{k+1}$ then there exists $s > t$ such that $u(s) \in A_k$. To this end, it is enough to show, that if $u(t) \in A_{k+1} \backslash A_k$, then there exists $s > t$ such that $u(s) \in A_k$.

**Step 3.1.** We first prove that if $u_j(t) < d_j + (k+1)\theta c_j$ for every $j \in \{1, \ldots, n\}$ and every $t \geqslant t_0$, then for every $i \in \{1, \ldots, n\}$ there exists $s > t_0$ such that $u_i(s) < d_i + k\theta c_i$. To this end suppose that $u_i(t) \geqslant d_i + k\theta c_i$ for all $t \geqslant t_0$. Then for every $t \geqslant t_0$

$$u_i(t) = u_i(t_0) \exp\left( \int_{t_0}^t a_i(r) - b_{ii}(r)u_i(r) - \sum_{j=1, j \neq i}^n b_{ij}(r)u_j(r)\, dr \right)$$

$$\leqslant u_i(t_0) \exp\left( \int_{t_0}^t a_i(r) - b_{ii}(r)d_i - b_{ii}(r)c_j k\theta - \sum_{j=1, j \neq i}^n b_{ij}(r)(d_j + (k+1)\theta c_j)\, dr \right)$$

$$= u_i(t_0) \exp\left( \int_{t_0}^t a_i(r) - b_{ii}(r)d_i - \sum_{j=1, j \neq i}^n b_{ij}(r)(d_j + \theta c_j) - k\theta \sum_{j=1}^n b_{ij}(r)c_j\, dr \right)$$

$$\leqslant u_i(t_0) \exp\left( \int_{t_0}^t a_i(r) - b_{ii}(r)d_i - \sum_{j=1, j \neq i}^n b_{ij}(r)(d_j + \theta c_j)\, dr \right)$$

If $i \in I$, then, by $(B)$, $u_i(t) \leqslant u_i(t_0) \exp(-\varepsilon(t - t_0))$, a contradiction with $u_i(t) \geqslant d_i + k\theta c_i$ for every $t \geqslant t_0$. If $i \in J$, then as $b_{ii}(r)d_i > 0$, it follows by $(B)$ that also $u_i(t) \leqslant u_i(t_0) \exp(-\varepsilon(t - t_0))$ and we arrive at the same contradiction.

**Step 3.2.** We now show that if for some $t$ and for some $i$ we have $u_i(t) < d_i + k\theta c_i$ and for every $s \geqslant t$ and every $j \in \{1, \ldots, n\}$ we have $u_j(s) < d_j + (k+1)\theta c_j$, then for every $s \geqslant t$ we also have $u_i(s) < d_i + k\theta c_i$. Assume the contrary, then there exists $t^* > t$ such that $u_i(t^*) = d_i + k\theta c_i$ and $\dot{u}_i(t^*) \geqslant 0$. This means that

$$a_i(t^*) \geqslant b_{ii}(t^*)u_i(t^*) + \sum_{j=1, j \neq i}^n b_{ij}(t^*)u_j(t^*)$$

$$\geqslant b_{ii}(t^*)d_i + b_{ii}(t^*)k\theta c_i + \sum_{j=1, j \neq i}^n b_{ij}(t^*)d_j + \sum_{j=1, j \neq i}^n b_{ij}(t^*)(k+1)\theta c_j$$

$$= b_{ii}(t^*)d_i + k\theta \sum_{i=1}^n b_{ii}(t^*)c_i + \sum_{j=1, j \neq i}^n b_{ij}(t^*)(d_j + \theta c_j)$$

$$\geqslant b_{ii}(t^*)d_i + \sum_{j=1, j \neq i}^n b_{ij}(t^*)(d_j + \theta c_j).$$

Now if $i \in I$ we immediately get the contradiction with $(B)$, and if $i \in J$ the last estimate implies that $a_i(t^*) \geqslant \sum_{j=1, j \neq i}^n b_{ij}(t^*)(d_j + \theta c_j)$ which also, by $(B)$ leads to a contradiction which ends the proof of Step 3.

**Step 4.** For every initial condition $u_0 = (u_1(t_0), u_2(t_0), \ldots, u_n(t_0)) > 0$ there exists some $k \geqslant 0$ such that $u_0 \in A_k$, so, repeating the process above, we know that there exists some $t^*$ such that $u(t^*) \in A_0$, and then

$$u_i(t) < d_i$$

for every $i = 1, \ldots, N$ and every $t \geqslant t^*$.                                    $\square$

**Remark 2.** *Redheffer in [Red96, Lemma 5] obtains a similar result. However, note that we do not need that*

$$\sum_{j=1}^n b_{ij}(t)d_j > 0 \quad for \ \ i \in J.$$

*The row diagonal dominance condition* $(H1)$ *holds with constants* $c_j$, *but it does not have to hold with the constants* $d_j$, *by which* $b_{ij}$ *and* $a_j$ *relate.*

**Lemma 4.2.** *Assume* $(H1)$ *and* $(B)$. *Let* $u = (u_1, u_2, ..., u_n)$ *be solution of* (LV-$n$) *such that* $u_i(t_0) > 0$ *for* $i = 1, .., n$. *If* $i \in J$ *then* $u_i \to 0$ *when* $t \to \infty$.

*Proof.* From Lemma 4.1 it follows that there exists $t^* \geqslant t_0$ such that for every $t \geqslant t^*$, we have $u_i(t) < d_i$ for $i = 1, ..., n$. Then, if only $k \in J$

$$\hat{u}_k(t) = \hat{u}_k(t^*) \exp\left(\int_{t^*}^t a_k(s) - \sum_{j=1}^n b_{kj}(s) u_j(s) ds\right) \leqslant d_k \exp\left(\int_{t^*}^t a_k(s) - \sum_{j=1, j \neq k}^n b_{kj}(s) d_j \, ds\right)$$

$$\leqslant d_k \exp\left(-\varepsilon(t - t^*)\right) \to 0,$$

whence $u_k(t) \to 0$ as $t \to \infty$. □

## 4.2 Backward unboundedness of solutions in case of extinction

In this subsection we prove that if the set $J$ is not empty, then any solution with all positive initial data must be backward unbounded.

**Lemma 4.3.** *Assume* $(H1)$ *and* $(B)$ *and let* $J$ *be nonempty. Every solution* $u$ *of* (LV-$n$) *with the initial condition* $u(t_0) > 0$ *is backward unbounded.*

*Proof.* We take a solution $u$ with the initial condition $u(t_0) \in A_0$, where $A_0$ is given by (10) with $k = 0$. We take $i \in J$, then

$$u_i' = u_i\left(a_i - b_{ii} u_i - \sum_{j \neq i} b_{ij} u_j\right) \leqslant u_i\left(a_i - \sum_{j \neq i} b_{ij} d_j\right) - b_{ii} u_i^2 \leqslant u_i\left(\sum_{j \neq i} b_{ij} c_j \theta - \varepsilon\right) - b_{ii} u_i^2 \leqslant -\varepsilon u_i.$$

Since $A_0$ is positively invariant, the straightforward application of the Gronwall lemma yields together with the fact that $u(t_0) \in A_0$

$$d_i \geqslant u_i(t_0) \geqslant u_i(t) e^{\varepsilon(t - t_0)}.$$

We deduce that there exists $t_0^*$ such that $u(t_0^*) \notin A_0$, otherwise we would get a contradiction by passing with $t_0 \to -\infty$.

Now we suppose that the solution has the initial condition $u(t_0) \in A_{k+1} \backslash A_k$. We are going to prove that there exists $t^* \leqslant t_0$ such that, $u(t^*) \in A_{k+2} \backslash A_{k+1}$, and then we would get the backward unboundedness of the solution by iteration. For contradiction assume that for every $t \leqslant t_0$ we have $u(t) \in A_{k+1}$. Then, as every set in family $\{A_k\}_{k \geqslant 0}$ is positively invariant, for every $t \leqslant t_0$ $u(t) \notin A_k$, and, arguing as in Step 3.2 of Lemma 4.1 there exists $i \in \{1, 2, ..., n\}$ such that $u_i(t) \geqslant d_i + k \theta c_i$ for every $t \leqslant t_0$.

Suppose that $i \in I$. We prove that it must be $\dot{u}_i(t) \leqslant -\varepsilon$ for every $t \leqslant t_0$. Indeed, if there exists $t_1 \leqslant t_0$ such that the opposite inequality holds, then at that time

$$a_i > b_{ii} u_i + \sum_{j \neq i} b_{ij} u_j - \varepsilon \geqslant b_{ii} d_i + b_{ii} k \theta c_i + \sum_{j \neq i} (b_{ij} d_j + b_{ij} k \theta c_j + b_{ij} \theta c_j) - \varepsilon \geqslant b_{ii} d_i + \sum_{j \neq i} b_{ij} (d_j + \theta c_j) - \varepsilon$$

and we get a contradiction with $(B)$. Then, there exists $t^* \leqslant t_0$ such that $u_i(t^*) \geqslant d_i + (k+2)\theta c_i$ so it can not be $u(t) \in A_{k+1}$ for all the time in the past, a contradiction.

Analogous argument with $i \in J$ yields the inequality

$$a_i > \sum_{j \neq i} b_{ij} d_j + b_{ij} \theta c_j - \varepsilon,$$

which also leads to a contradicion. $\qquad \square$

## 4.3 The case of extinction of all species except one

We have given conditions which guarantee the extinction of given a subset of species. In this section we study the asymptotic behavior of the persistent species indexed by $I$. We prove that, as time goes to infinity, the quantities of these species tend to the unique separated from zero and infinity solutions of the subsystem given by the equations with index $I$. We start from the situation when only one species persists. Hence, we provide the conditions over the vector $a(t)$ so that the trajectory $(u_1^*, 0, ..., 0)$ is the global attractive solution of the system (LV-$n$). The function $u_1^*$ is the solution of the logistic equation

$$\dot{u}_1(t) = u_1(t)(a_1(t) - b_{11}(t)u_1(t))$$

separated from zero and infinity, which, by Theorem 2.1 is unique under assumptions $(B1)$ below (which implies $(B)$ for one dimensional system) and $(H2)$. Hence, we need to impose $b_{11} > 0$ (already imposed by the initial assumptions) and the version of $(B)$ with $I = \{1\}$ and $J = \{2, ..., n\}$, namely we assume that there exists a number $\bar{d}_1 > 0$, a vector $d > 0$, and two numbers $\varepsilon, \theta > 0$ such that,

$$\begin{cases} b_{11}\bar{d}_1 + \varepsilon \leqslant a_1 \leqslant b_{11}d_1 + \sum_{j=2}^n b_{1j}(d_j + c_j\theta) - \varepsilon \\ a_i \leqslant \sum_{j=1, j \neq i}^n b_{ij}(d_j + c_j\theta) - \varepsilon \text{ for } i = 2, ..., n \end{cases} \tag{B1}$$

We skip the proof of the following result which exactly follows the lines of the proof of [Ahm93, Lemma 4], where, however, only two species, one that persists and one that decays, are considered.

**Lemma 4.4.** *Assume $(B1)$. Let be $(u_1, u_2, ..., u_n)$ a solution of (LV-$n$) such that $\bar{d}_1 < u_1(t_0) < d_1$ and $0 < u_i(t_0) < d_i$ for every $i = 2, ..., n$, then $u_1^* - u_1 \to 0$ as time tends to infinity.*

Finally, combining the above lemma with Lemma 4.1 and Lemma 4.2 we have the following result on the asymptotic behavior.

**Theorem 4.5.** *Suppose the hypothesis $(H1)$ and $(B1)$ and let $|b_{1j}|^U < \infty$ for $j = \{2, ..., n\}$, if $(u_1, u_2, ..., u_n)$ is a solution with the initial conditions $u_i(t_0) > 0$ for $i = 1, ..., n$, then $u_i \to 0$ for $i = 2, ..., n$ and $u_1 - u_1^* \to 0$ as time tends to infinity.*

*Proof.* By Lemma 4.1, we can choose $t^*$ as the initial time, where $\bar{d}_1 < u_1(t^*) < d_1$ and $0 < u_i(t^*) < d_i$ for $i = 2, ..., n$. Then, by Lemma 4.4, we have that $u_1(t) - u_1^*(t) \to 0$ and by Lemma 4.2 we deduce $u_i(t) \to 0$ when $t \to \infty$ for $i = 2, ..., n$. $\qquad \square$

### 4.4 The case of extinction of one species

We continue with the analysis of the case when only one species goes to extinction. Hence, we provide the conditions on the coefficients $a_i(t), b_{ij}(t)$ which guarantee that the trajectory $(u_1^*, u_2^*, ..., u_{n-1}^*, 0)$ is the global attractive solution when $t \to \infty$ for (LV-$n$), and where $(u_1^*, u_2^*, ..., u_{n-1}^*)$ is the solution of the Lotka–Volterra $(n-1)$-dimensional system obtained by removing the last variable.

Now need to impose both $(H1)$, and the version of $(B)$ with $I = \{1, ..., n-1\}$ and $J = \{n\}$, namely the existence of vectors $d, \bar{d} \in \mathbb{R}^n$ and parameters $\theta, \varepsilon > 0$ such that

$$\begin{cases} b_{ii}\bar{d}_i + \varepsilon \leqslant a_i \leqslant b_{ii}d_i + \sum_{j=1; j \neq i}^{n} b_{ij}(d_j + c_j\theta) - \varepsilon & \text{for} \quad i = 1, ..., n-1 \\ a_n \leqslant \sum_{j=1}^{n-1} b_{nj}(d_j + c_j\theta) - \varepsilon \end{cases} \quad (B2)$$

for every $t \in \mathbb{R}$. Also it has to hold the column diagonal dominance for the first $n-1$ coordinates, so we are going to say that holds $(H2)_{n-1}$.

**Lemma 4.6.** *Assume $(H2)_{n-1}$ and $(B2)$. Let be $\bar{u}$ be a solution of an $(n-1)$-dimensional Lotka–Volterra system obtained by setting the last variable to zero $(u_n = 0)$, with the initial data $\bar{u}(t_0) > 0$. There exists a unique solution $u^*$ of $(n-1)$-dimensional system, which is separated from zero and infinity and such that*

$$|\bar{u}(t) - u^*(t)| \to 0 \quad when \quad t \to \infty$$

*and*

$$\bar{d}_i \leqslant u_i^*(t) \leqslant d_i$$

*for $i = 1, ..., n-1$, and for all $t \in \mathbb{R}$.*

*Proof.* We use Theorem 2.1 for dimension $n-1$. $\qquad \square$

**Lemma 4.7.** *Assume $(H1)$, $(H2)_{n-1}$ and $(B2)$ and let $|b_{in}|^U < \infty$ for $i \in \{1, ..., n-1\}$. Let $u^* = (u_1^*, u_2^*, ..., u_{n-1}^*)$ be the solution of $(n-1)$-dimensional Lotka–Volterra system given by Lemma 4.6 and let $\hat{u} = (u_1^*, u_2^*, ..., u_{n-1}^*, 0)$. There exists a constant $\delta > 0$ such that if for some $t_0 \in \mathbb{R}$ we have $|\hat{u}(t_0) - u_0| < \delta$, then $\lim_{t \to \infty} |\hat{u}(t) - u(t)| = 0$ where $u$ solves (LV-$n$) with $u(t_0) = u_0$.*

*Proof.* We start by studying the system linearized around $\hat{u} = (u_1^*, ..., u_{n-1}^*, 0)$. We write

$$w(t) = u(t) - \hat{u}(t)$$

where $u(t)$ is a solution with the initial condition in a neighbourhood of $\hat{u}$. Then we obtain

$$w'(t) = \mathcal{M}(t)w(t) + R(w(t), t),$$

where $\mathcal{M}(t)$ and $R(w, t)$ are as in (8). The linearized system has the form (9) with $C(t)$ being an $(n-1) \times 1$ which is bounded from the bound on $b_{in}$. Moreover $B(t)$ is $1 \times 1$ matrix with its entry given by $a_n(t) - \sum_{j=1}^{n-1} b_{n,j}(t)u_j^*(t)$. We observe that the one dimensional system $v_n' = B(t)v_n$ has an exponential dichotomy with projection $P(t) = 1$ since by $(B2)$, it holds $B(t) \leqslant -\varepsilon < 0$. Moreover, by Lemma 3.3, the system $(v_1, ..., v_{n-1})' = A(t)(v_1, ..., v_{n-1})$ has an exponential dichotomy with $P(t) = I_{(n-1) \times (n-1)}$, the $(n-1) \times (n-1)$ identity.

Hence, by Corollary 3.2 and results of [BP15] recalled in Section 3 the linearized system (9) admits an exponential dichotomy with the projection given by $P = I_{n \times n}$. Is easy to see that the time dependent projection $P(t)$ is given by $P(t) = M(t, 0)PM(0, t) = I_{n \times n}$. So

$$|v(t)| \leqslant K|v(s)|e^{-\gamma(t-s)} \quad \text{for all} \quad t, s \in \mathbb{R},$$

for some constants $K, \gamma > 0$, and by [Cop65, page 70] we get the local asymptotic stability of $\hat{u}$ given in the assertion of the lemma.

<div style="text-align: right">□</div>

In the next result we establish the global asymptotic stability. We skip the proof because it exactly follows the lines of [AL98, Theorem 2.3].

**Theorem 4.8.** *Assume $(H1)$, $(H2)$ and $(B2)$ and let $|b_{ij}|^U < \infty$ and $|a_i|^U < \infty$ for $i, j \in \{1, \dots, n\}$. Let be $u$ a solution of $(\mathrm{LV}\text{-}n)$ such that $u(t_0) > 0$. Then*

$$\lim_{t \to \infty} |u(t) - \hat{u}(t)| = 0$$

*where $\hat{u} = (u_1^*, \dots, u_{n-1}^*, 0)$, $(u_1^*, \dots, u_{n-1}^*)$ being the globally asymptotically stable positive solution solution of $(n-1)$-dimensional system given in Lemma 4.6.*

# 5 Structure of the attractor for non-autonomous logistic equation

We start from the analysis for the non-autonomous problem in one dimension. Although the results of this section are mostly known, we need them for the study of problems in two and three dimensions. Thus, the aim of this chapter is the study of

$$u'(t) = u(t)(a(t) - b(t)u(t)), \tag{11}$$

where $a, b \in C(\mathbb{R}), b > 0$.

## 5.1 Permanence

We need the assumption

$$\overline{d}b(t) \leqslant a(t) \leqslant db(t) \quad \text{and} \quad b^L > 0, \tag{12}$$

with constants $\overline{d}, d > 0$. The next result follows from [Red96, Theorem 1. (ii)], cited above as Theorem 2.1.

**Lemma 5.1.** *The function $u(t) = 0$ for $t \in \mathbb{R}$ is a solution of (11). Moreover, there exists a function $\overline{d} \leqslant u^*(t) \leqslant d$ for $t \in \mathbb{R}$ which is a solution to (11). This is the unique complete solution separated away from zero and infinity.*

Next lemma provides the characterization of the asymptotic behavior of solutions to (11) other than $u^*$.

**Lemma 5.2.** *If $u : \mathbb{R} \to \mathbb{R}$ is a solution to (11) with $u(t_0) \geqslant 0$ then exactly one of the four possibilities below holds:*

(a) $u(t) = 0$ for every $t \in \mathbb{R}$,

(b) $u = u^*$,

(c) if $u(t_0) \in (0, u^*(t_0))$ then $\lim_{t \to -\infty} u(t) = 0$ and $\lim_{t \to \infty} u^*(t) - u(t) = 0$,

(d) if $u(t_0) > u^*(t_0)$ then $\lim_{t \to -\infty} u(t) = \infty$ and $\lim_{t \to \infty} u(t) - u^*(t) = 0$.

*Proof.* If we take the initial data $u(t_0) \in (0, u^*(t_0))$ then it is clear that $u(t) \in (0, u^*(t))$ for every $t$. If $\inf_{t \in \mathbb{R}} u(t) > 0$, then by [Red96, Theorem 1 (v)] we get a contradiction because the solution bounded away from zero must be unique. So we suppose that $\inf_{t \in \mathbb{R}} u(t) = 0$. If there exists a time $t_1 \in \mathbb{R}$ such that $u(t_1) = 0$, then $u = 0$ by uniqueness of solution. Then there exists a decreasing sequence $\{t_n\}$ such that $t_n \to -\infty$, $u(t_n) \to 0$ and $u(t_n) < \overline{d}$. If for some $t \in (t_{n+1}, t_n)$, we have $u(t) > u(t_n)$ then there exists $t^* \in [t, t_n)$ such that $u(t^*) > u(t_n)$ and on $(t^*, t_n)$ the function $u$ is strictly less that $\overline{d}$. Hence

$$0 > u(t_n) - u(t^*) = \int_{t^*}^{t_n} u(s)(a(s) - b(s)u(s)) \, ds \geqslant \int_{t^*}^{t_n} u(s)(a(s) - b(s)\overline{d}) ds \geqslant 0,$$

a contradiction, and hence $\lim_{t \to -\infty} u(t) = 0$. Convergence to $u^*$ as $t \to \infty$ follows from [Red96, Theorem 1 (vi)]. If $u(t_0) > u^*(t_0)$ then analogously as for the case (c) we have the convergence at $+\infty$ and the existence of a decreasing sequence $t_n \to -\infty$ such that $u(t_n) \to \infty$ and $u(t_n) > d$. If $u(t) < u(t_n)$ for $t \in (t_{n+1}, t_n)$ then there exists $t^* \in [t, t_n)$ such that $u(t^*) < u(t_n)$ and $u$ is strictly greater than $d$ on $[t^*, t_n)$ hence, analogously as in the case (c), $\lim_{t \to -\infty} u(t) = \infty$. □

We finish this section by recalling two results on one dimensional problem useful in further analysis. We first recall a result that the linearized problem has the exponential dichotomy with its only direction being exponentially stable.

**Lemma 5.3.** *There exist $\delta > 0$ and $M > 0$ such that if $w : \mathbb{R} \to \mathbb{R}$ solves the following one-dimensional problem linearized around $u^*$*

$$w'(t) = w(t)(a(t) - 2b(t)u^*(t)),$$

*then we have*

$$|w(t)| \leqslant M|w(t_0)|e^{-\delta(t-t_0)} \quad \text{for every} \quad t \geqslant t_0.$$

*Proof.* The result is a direct application of Lemma 3.3 which in turn follows from [AL98, Lemma 3.6]. □

Finally we establish the invariance and monotonicity result for the solution of one-dimensional problem (11).

**Lemma 5.4.** *If $u(t_0) \geqslant \overline{d}$ then $u(t) \geqslant \overline{d}$ for every $t \geqslant t_0$ and if $u(t_0) \leqslant d$ then $u(t) \leqslant d$ for every $t \geqslant t_0$. If $0 < u(t_0) < \overline{d}$ then $u(t) \geqslant u(t_0)$ for every $t \geqslant t_0$ and if $u(t_0) > d$ then $u(t) \leqslant u(t_0)$ for every $t \geqslant t_0$.*

*Proof.* The first part follows from [Red96, Theorem 1. (i)]. For the second part it is enough to prove that if $0 < u(t) < \overline{d}$ then $u'(t) > 0$. Indeed if $u(t_0) \in (0, \overline{d})$

$$u'(t) = u(t)\left(a(t) - b(t)u(t)\right) > u(t)\left(a(t) - b(t)\left(\frac{a}{b}\right)^L\right) = u(t)b(t)\left(\frac{a(t)}{b(t)} - \left(\frac{a}{b}\right)^L\right) \geqslant 0.$$

In turn, if $u(t) > d$, then

$$u'(t) = u(t) \left(a(t) - b(t)u(t)\right) < u(t)\left(a(t) - b(t)\left(\frac{a}{b}\right)^U\right) = u(t)b(t)\left(\frac{a(t)}{b(t)} - \left(\frac{a}{b}\right)^U\right) \leqslant 0,$$

and the proof is complete. □

## 5.2 Extinction

To get the criterion on the extinction of the single species we give a theorem that is a direct application of Lemma 4.2.

**Lemma 5.5.** *If $a^U < 0$ then for every solution of* (11) *with $u(t_0) > 0$ we have $\lim_{t \to \infty} u(t) = 0$ and $\lim_{t \to -\infty} u(t) = \infty$.*

*Proof.* To get the assertion $\lim_{t \to -\infty} u(t) = \infty$ observe that since

$$\frac{u'(t)}{u(t)} = a(t) - b(t)u(t),$$

we deduce that

$$\ln(u(t)) - \ln(u(t_0)) = a(t) - b(t)u(t) \leqslant a^U.$$

Hence

$$u(t_0) \geqslant u(t)e^{-a^U(t-t_0)},$$

and the proof is complete. □

# 6 Attractor for non-autonomous Lotka-Volterra 2-D system

In Section 2 we recalled the result of [Red96] on the existence of complete trajectories which are bounded and separated from zero. In particular we have shown condition which guarantees the existence of such solution. This condition can be used for the whole $n$-dimensional system or for its subsystems obtained by setting some of the variables $u_i$ to zero. Using this theorem, combined with the results of Section 4, we present the conditions under which one can characterize the structure of the non-autonomous attractor for cooperative Lotka–Volterra problem in two dimensions. We will consider three cases depending on the globally asymptotically stable solution: either the one with both nonzero components, or the one with one nonzero component, or the one with both zeros. Hence we will consider the system

$$\begin{cases} u'_1 = u_1(a_1(t) - b_{11}(t)u_1 - b_{12}(t)u_2) \\ u'_2 = u_2(a_2(t) - b_{21}(t)u_1 - b_{22}(t)u_2), \end{cases} \tag{LV-2}$$

with $b_{11}^L, b_{22}^L > 0$ and $b_{12}, b_{21} \leqslant 0$.

We will use the notation $S(t, \tau)u_0$ for $\tau \in \mathbb{R}$ and $t \geqslant \tau$ to denote the process associated to (LV-2), i.e. the family of maps which assign to the initial data $u_0$ at time $\tau$ the solution of (LV-2) at time $t$.

## 6.1 Structure of attractor for the case of permanence

We will first study the case when there exists the complete solution bounded away from zero and infinity on both variables, which attracts all solutions with nonzero initial data as $t \to \infty$. In such case, of the permanence of the two species, the structure of non-autonomous attractor is depicted in Figure 1.

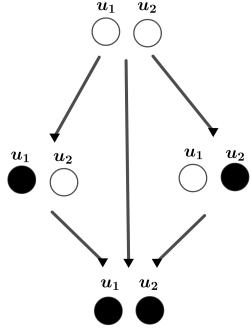

Figure 1: Four complete solutions and their connections for the two dimensional problem with globally asymptotically stable state with coexistence of both species. Black dot means the strictly positive function while white dot relates with the function identically equal to zero.

Following Theorem 2.1 we have to impose the following conditions

$$\left\{ \begin{array}{l} (\overline{c}_1 b_{11} + \overline{c}_2 b_{21})^L > 0, \\ (\overline{c}_2 b_{22} + \overline{c}_1 b_{12})^L > 0, \end{array} \right. \quad \left\{ \begin{array}{l} \overline{d}_1 b_{11} \leqslant a_1 \leqslant d_1 b_{11} + d_2 b_{12}, \\ \overline{d}_2 b_{22} \leqslant a_2 \leqslant d_1 b_{21} + d_2 b_{22}, \end{array} \right. \tag{13}$$

with some constants $\overline{d}_1, \overline{d}_2, d_1, d_2, \overline{c}_1, \overline{c}_2 > 0$. These bounds are exactly the conditions $(A)$ and $(H2)$ for the two-dimensional case. We also assume that

$$|b_{12}|^U, |b_{21}|^U, a_1^U, a_2^U < \infty. \tag{14}$$

Note that the above assumptions in particular imply that $a_1, a_2 > 0$. The next result follows from [Red96, Theorem 1. (ii)], cited above as Theorem 2.1.

**Lemma 6.1.** *There exists a function $u^* = (u_1^*, u_2^*)$ defined for $t \in \mathbb{R}$, the complete solution of* (LV-2) *such that $\overline{d}_i \leqslant u_i \leqslant d_i$ for $i = 1, 2$. This is the unique complete trajectory bounded away from zero and infinity in both variables. Moreover, there exist functions $(\hat{u}_1, 0)$, $(0, \hat{u}_2)$ and $u(t) = (0, 0)$ defined for $t \in \mathbb{R}$ that are complete solutions to* (LV-2), *such that $\overline{d}_i \leqslant \hat{u}_i \leqslant d_i$ for $i = 1, 2$.*

The next Lemma follows directly from Lemma 5.2.

**Lemma 6.2.** *If the function $u : \mathbb{R} \to \mathbb{R}^2$ is a solution to* (LV-2) *with $u_1(t_0) \geqslant 0$ and $u_2(t_0) = 0$, then exactly one of the four possibilities below holds:*

(a) $u(t) = (0,0)$ *for every* $t \in \mathbb{R}$,

(b) $u = (\hat{u}_1, 0)$,

(c) *if* $u_1(t_0) \in (0, \hat{u}_1(t_0))$ *then* $\lim_{t \to -\infty} u_1(t) = 0$, $\lim_{t \to \infty} \hat{u}_1(t) - u_1(t) = 0$ *and* $u_2(t) = 0$ *for* $t \in \mathbb{R}$,

(d) *if* $u_1(t_0) > \hat{u}_1(t_0)$ *then* $\lim_{t \to -\infty} u_1(t) = \infty$, $\lim_{t \to \infty} u_1(t) - \hat{u}_1(t) = 0$, *and* $u_2(t) = 0$ *for* $t \in \mathbb{R}$.

The result analogous to Lemma 6.2 holds, and allows us to establish the connection from $(0,0)$ to $(0, \hat{u}_2)$. In order to obtain the structure as in Fig. 1, it remains to establish connections from $(0, \hat{u}_2)$, $(\hat{u}_1, 0)$ and $(0,0)$ to $(u_1^*, u_2^*)$.

**Theorem 6.3.** *There exists the trajectory of* (LV-2) *denoted by* $z = (z_1, z_2) : \mathbb{R} \to \mathbb{R}^2$ *such that*

$$\lim_{s \to -\infty} |(z_1(s), z_2(s)) - (\hat{u}_1(s), 0)| = 0.$$

*and*

$$\lim_{s \to \infty} |(z_1(s), z_2(s)) - (u_1^*(s), u_2^*(s))| = 0.$$

*Analogous result holds for* $(0, \hat{u}_2)$.

*Proof.* We prove that the solution $(\hat{u}_1, 0)$ is locally unstable, i.e. its non-autonomous unstable manifold

$$W^u((\hat{u}_1, 0)) = \{(t, (w_1, w_2)) : \text{ there exists a solution } z : \mathbb{R} \to \mathbb{R}^2 \text{ such that}$$
$$z(t) = (w_1, w_2) \text{ and } \lim_{s \to -\infty} |z(s) - (\hat{u}_1(s), 0)| = 0\},$$

is nonempty and intersects the interior of the positive quadrant. As we have seen in Section 2, cf Theorem 2.1, the permanent solution $(u_1^*(s), u_2^*(s))$ attracts forward in time the solution with the initial data in this intersection. We start by proving the existence of the local unstable manifold. To this end we first study the system linearized around $(\hat{u}_1, 0)$. We write

$$w(t) = u(t) - \hat{u}(t)$$

where $\hat{u}(t) = (\hat{u}_1(t), 0)$, and $u(t)$ is a solution with initial condition in a neighbourhood of $\hat{u}$. Then we obtain

$$w'(t) = \mathcal{M}(t)w(t) + R(w(t), t)$$

where $\mathcal{M}(t)$ is the derivative of the vector field at $\hat{u}(t)$ and $R(t)$ is a remainder which is of higher (quadratic) order with respect to $w$. In the 2D case we have the following ODE

$$w'(t) = \begin{pmatrix} a_1(t) - 2b_{11}(t)\hat{u}_1(t) & -b_{12}(t)\hat{u}_1(t) \\ 0 & a_2(t) - b_{21}(t)\hat{u}_1(t) \end{pmatrix} w(t) + \begin{pmatrix} -b_{11}(t)w_1(t)^2 - b_{12}(t)w_1(t)w_2(t) \\ -b_{21}(t)w_1(t)w_2(t) - b_{22}(t)w_2(t)^2 \end{pmatrix}$$
$$(15)$$

The linearized system has the form

$$\begin{cases} v'(t) = \begin{pmatrix} a_1(t) - 2b_{11}(t)\hat{u}_1(t) & -b_{12}(t)\hat{u}_1(t) \\ 0 & a_2(t) - b_{21}(t)\hat{u}_1(t) \end{pmatrix} v(t) = \begin{pmatrix} A(t) & C(t) \\ 0 & B(t) \end{pmatrix} v(t). \end{cases} \qquad (16)$$

where $A(t) = a_1(t) - 2b_{11}(t)\widehat{u}_1(t)$, $C(t) = -b_{12}(t)\widehat{u}_1(t)$ and $B(t) = a_2(t) - b_{21}(t)\widehat{u}_1(t)$. We observe that the one dimensional system $v_2' = B(t)v_2$ has an exponential dichotomy with projection $P_B(t) \equiv 0$ since by (13) it holds $B(t) \geqslant \delta > 0$ and for every $t \geqslant s$

$$v_2(t) = v_2(t)e^{\int_s^t a_2(r) - b_{21}(r)\widehat{u}_1(r)ds} \Rightarrow |v_2(t)| \geqslant |v_2(s)|e^{\delta(t-s)}.$$

We also observe that by Lemma 5.3 the system $v_1' = A(t)v_1$ has an exponential dichotomy with $P_A(t) = I$. Then, by results recalled in Section 3, since by (14) and Lemma 6.1 the function $C(t)$ is bounded, the system (16) admits an exponential dichotomy. According to (7), the projections are given by

$$P^+ = \begin{pmatrix} 1 & 0 \\ 0 & 0 \end{pmatrix} \text{ on } \mathbb{R}^+ \text{ and } P^- = \begin{pmatrix} 1 & L^- \\ 0 & 0 \end{pmatrix} \text{ on } \mathbb{R}^-.$$

The projection

$$P = P^- = \begin{pmatrix} 1 & L^- \\ 0 & 0 \end{pmatrix}$$

has the same range as $P^+$ and hence the system (16) has exponential dichotomy with $P$. For the time $t \in \mathbb{R}$ the associated projection is given by $P(t) = M(t,0)PM(0,t)$, where $M(t,\tau)$ is the fundamental matrix of (16), which yields, as both $M(t,0)$ and $M(0,t)$ are upper triangular and inverse to each other,

$$P(t) = \begin{pmatrix} 1 & L(t) \\ 0 & 0 \end{pmatrix},$$

for certain bounded function $L(t)$. In particular $\dim \ker P(t) = 1$ and $\dim \text{range } P(t) = 1$. Moreover,

$$\text{range}(P(t)) = \left\{ \beta \begin{pmatrix} 1 \\ 0 \end{pmatrix} : \beta \in \mathbb{R} \right\} \text{ and range}(I - P(t)) = \left\{ \alpha \begin{pmatrix} -L(t) \\ 1 \end{pmatrix} : \alpha \in \mathbb{R} \right\}.$$

We will use [KR11, Theorem 6.10 and Exercise 6.11]. There exists a neighborhood $U$ of zero in $\mathbb{R}^2$ and a continuous function $\Sigma^- : \mathbb{R} \times U \to \mathbb{R}^2$, a local non-autonomous unstable manifold, such that $\Sigma^-(t,x) = \Sigma^-(t,(I - P(t))x) \in \text{range}(P(t))$, $\Sigma^-(t,0) = 0$, $\lim_{|x| \to 0} \frac{\Sigma^-(t,x)}{|x|} = 0$. Moreover, if for some $t \in \mathbb{R}$ and $x \in \text{range}(I - P(t)) \cap U$ we have $y = x + \Sigma^-(t,x)$, then, provided $|y|$ is small enough, the backward trajectory $w : (-\infty, t] \to \mathbb{R}^2$ of (15) with $w(t) = y$ belongs to the graph of $\Sigma^-$, i.e. $\Sigma^-(\tau, (I - P(\tau))w(\tau)) = P(\tau)w(\tau)$ and satisfies $|w(\tau)| \leqslant Ce^{\delta(\tau - t)}$ for every $\tau \leqslant t$. The points in the graph of $\Sigma^-$ have the form

$$y = x + \Sigma^-(t,x) = \alpha \begin{pmatrix} -L(t) \\ 1 \end{pmatrix} + \Sigma^- \left( t, \alpha \begin{pmatrix} -L(t) \\ 1 \end{pmatrix} \right) = \begin{pmatrix} \Sigma^- \left( t, \alpha \begin{pmatrix} -L(t) \\ 1 \end{pmatrix} \right) - \alpha L(t) \\ \alpha \end{pmatrix}. \quad (17)$$

Because $L(t)$ is bounded, there exists a neighborhood of zero in $\mathbb{R}$ such that for $\alpha$ in this neighborhood, we have $x = \alpha \begin{pmatrix} -L(t) \\ 1 \end{pmatrix} \in U$. We take small positive $\alpha$. The second coordinate of $y$ is equal to $\alpha$ and hence positive, while the first coordinate can be made arbitrarily small, so that after adding to the solution $(\widehat{u}_1(t), 0)$, with $\widehat{u}_1$ separated away from zero, the first coordinate of the sum is also positive. Hence, for any $t \in \mathbb{R}$ there exists the point with second coordinate positive and first arbitrarily small, at time $t$, which is backward exponentially attracted to zero, and in terms of the original system, there exists the point at any $t \in \mathbb{R}$ with both coordinates positive which is backwards exponentially attracted to $(\widehat{u}_1(t), 0)$. In view of [Red96, Theorem 1 (iii) and (vi)] this gives us the assertion of the Theorem. □

**Theorem 6.4.** *There exists a trajectory of* (LV-2) *denoted by* $y = (y_1, y_2) : \mathbb{R} \to \mathbb{R}^2$ *such that*

$$\lim_{s \to -\infty} |(y_1(s), y_2(s))| = 0.$$

*and*

$$\lim_{s \to \infty} |(y_1(s), y_2(s)) - (u_1^*(s), u_2^*(s))| = 0.$$

*Proof.* The linearized system at $(0, 0)$ has the form

$$\begin{cases} v'(t) = \begin{pmatrix} a_1(t) & 0 \\ 0 & a_2(t) \end{pmatrix} v(t), \end{cases} \tag{18}$$

which has the exponential dichotomy with $P(t) \equiv 0$. Thus, by [KR11, Theorem 6.10 and Exercise 6.11] there exists a sufficiently small neighborhood of $(0, 0)$ such that the constant function $\Sigma^- : \mathbb{R} \times \mathbb{R}^2 \to \{0\}$ is the local unstable manifold of $(0, 0)$ and thus every sufficiently small positive initial condition is exponentially backwards attracted to zero. Since it must be forwards attracted to $(u_1^*(s), u_2^*(s))$, the proof is complete. □

**Theorem 6.5.** *If a trajectory with a strictly positive initial consition is bounded both in the past and in the future then it must be one of trajectories described above.*

*Proof.* From [Red96, Theorem 1. (v)] we know that $u^*$ is the only solution on $\mathbb{R}$ that is bounded away from 0 and $\infty$. Thus, if the complete trajectory $u = (u_1, u_2)$ is bounded, not coincides with $u^*$, and none of the two coordinates is zero, then for at least one of its coordinates, say $u_1$, there must exist a decreasing sequence $t_n \to -\infty$ such that $u_1(t_n) \to 0$ and $u_1(t_n) < \overline{d}_1$. We shall prove that $\lim_{t \to -\infty} u_1(t) = 0$. We proceed analogously to the Lemma 5.2. If for some $t \in (t_{n+1}, t_n)$, we have $u_1(t) > u_1(t_n)$ then there exists $t^* \in [t, t_n)$ such that $u_1(t^*) > u_1(t_n)$ and on $(t^*, t_n)$ the function $u_1$ is strictly less that $\overline{d}_1$. Hence

$$0 > u_1(t_n) - u_1(t^*) = \int_{t^*}^{t_n} u_1(s)(a_1(s) - b_{11}(s)u_1(s) - b_{12}(s)u_2(s)) \, ds$$

$$\geqslant \int_{t^*}^{t_n} u_1(s)(a_1(s) - b_{11}(s)\overline{d}_1) ds \geqslant 0,$$

and we have the contradiction. Now we need to prove that if both $u_1$ and $u_2$ do not converge backwards to zero, say $u_1(t) \to 0$ as $t \to -\infty$ and $u_2(t)$ is separated from zero then it must be

$$\lim_{t \to -\infty} |u_2(t) - \widehat{u}_2(t)| = 0.$$

To this end assume that for some sequence $t_n \to -\infty$

$$\lim_{n \to \infty} |u_2(t_n) - \widehat{u}_2(t_n)| = c > 0.$$

For $t \geqslant 0$ we have

$$|u_2(t_n) - \widehat{u}_2(t_n)| \leqslant |S(t_n, t_n - t)(u_1(t_n - t), u_2(t_n - t)) - S(t_n, t_n - t)(0, \widehat{u}_2(t_n - t))|.$$

Since $u_2(t)$ is separated from zero, denote $K = (u_2)^L$ and $M = (u_2)^U$. For every $x \in [K, M]$

$$\begin{aligned}
|u_2(t_n) - \widehat{u}_2(t_n)| &\leqslant |S(t_n, t_n - t)(u_1(t_n - t), u_2(t_n - t)) - S(t_n, t_n - t)(0, x)| \\
&\quad + |S(t_n, t_n - t)(0, x) - S(t_n, t_n - t)(0, \widehat{u}_2(t_n - t))| \\
&\leqslant e^{\kappa t}|(u_1(t_n - t), u_2(t_n - t)) - (0, x)| + \frac{2 \max\{M, d_2\}^2}{\min\{K, \overline{d}_2\}} e^{-\delta \min\{K, \overline{d}_2\}t},
\end{aligned}$$

where we have used Lemma 2.2 and 2.3, and $\kappa$ depends on the bounds on $u_1, u_2$ and the coefficients of the problem. Pick $\varepsilon > 0$. We can find $t > 0$ such that $\frac{2 \max\{M, d_2\}^2}{\min\{K, \overline{d}_2\}} e^{-\delta \min\{K, \overline{d}_2\}t} < \frac{\varepsilon}{2}$. The sequence $u_2(t_n - t)$ has a convergent subsequence, let $x$ be a limit of this subsequence and pick $n$ large enough such that $Ce^{\kappa t}|(u_1(t_n - t), u_2(t_n - t)) - (0, x)| < \frac{\varepsilon}{2}$. Hence for $n$ large enough, on a subsequence,

$$|u_2(t_n) - \widehat{u}_2(t_n)| < \varepsilon,$$

a contradiction. $\qquad \square$

To finish this section we summarize the results obtained above, and we construct the full image of the atractor depicted in Figure 1.

**Theorem 6.6.** *Assuming* (13) *the system* (LV-2) *has the following trajectories* $u : \mathbb{R} \to \mathbb{R}^2$ *bounded both in the past and in the future:*

(a) $u(t) = (0, 0)$ *for* $t \in \mathbb{R}$,

(b) $u(t) = (\widehat{u}_1(t), 0)$ *and* $u(t) = (0, \widehat{u}_2(t))$, *corresponding to the unique solutions for one-dimensional subproblems bounded away from zero and infinity, given in Lemma 6.1.*

(c) *Solutions of type* $u(t) = (u_1(t), 0)$ *with initial condition* $0 < u_1(t_0) < \widehat{u}_1(t_0)$, *where* $\lim_{t \to -\infty} u_1(t) = 0$ *and* $\lim_{t \to \infty}(u_1(t) - \widehat{u}_1(t)) = 0$, *given in Lemma 6.2. Analogously,* $u(t) = (0, u_2(t))$ *with initial condition* $0 < u_2(t_0) < \widehat{u}_2(t_0)$, *where* $\lim_{t \to -\infty} u_2(t) = 0$ *and* $\lim_{t \to \infty}(u_2(t) - \widehat{u}_2(t)) = 0$.

(d) $u(t) = (u_1^*(t), u_2^*(t))$ *the unique solution with both nonzero coordinates bounded away from zero and infinity given in Lemma 6.1.*

(e) *Solutions of type* $u(t) = (u_1(t), u_2(t))$ *such that* $\lim_{t \to -\infty}(u_1(t), u_2(t)) = (0, 0)$ *and* $\lim_{t \to \infty}(u_1(t) - u_1^*(t), u_2(t) - u_2^*(t))$, *given in Theorem 6.4.*

(f) *Solutions of type* $u(t) = (u_1(t), u_2(t))$ *such that* $\lim_{t \to -\infty}(u_1(t), u_2(t)) = (\widehat{u}_1(t), 0)$ *and* $\lim_{t \to \infty}(u_1(t) - u_1^*(t), u_2(t) - u_2^*(t))$, *given in Theorem 6.3. Analogously,* $u(t) = (u_1(t), u_2(t))$ *such that* $\lim_{t \to -\infty}(u_1(t), u_2(t)) = (0, \widehat{u}_2(t))$ *and* $\lim_{t \to \infty}(u_1(t) - u_1^*(t), u_2(t) - u_2^*(t)) = 0$.

*Moreover any solution of* (LV-2) *which is bounded both in the past and in the future is one of the solutions described in items (a)-(f).*

## 6.2 Structure of attractor for the case of extinction of one species

Now we study is the extinction of one species. Hence, we consider (LV-2) with assumptions

$$
\left\{
\begin{array}{l}
(c_1 b_{11} + c_2 b_{12})^L > 0, \\
(c_1 b_{21} + c_2 b_{22})^L > 0,
\end{array}
\right.
\qquad
\left\{
\begin{array}{l}
b_{11}\overline{d}_1 + \varepsilon \leqslant a_1 \leqslant b_{11} d_1 + b_{12}(d_2 + \theta c_2) - \varepsilon, \\
a_2 \leqslant b_{21}(d_1 + \theta c_1) - \varepsilon,
\end{array}
\right.
\tag{19}
$$

with some constants $d_1, d_2, c_1, c_2, \overline{d}_1 > 0$ and constants $\varepsilon, \theta > 0$. Now, the first inequalities constitute the condition $(H1)$, while the second one in $(B1)$ for the two dimensional case. We also assume that $|b_{12}|^U < \infty$.

We prove that the dynamics of the problem is described by the diagram presented in Fig. 2. Indeed, it is sufficient to use Lemma 4.3, Theorem 4.5, Lemma 5.1, Lemma 5.2, and Lemma 5.5 to state the following theorem.

**Theorem 6.7.** *Assuming* (19) *the system* (LV-2) *has the following trajectories* $u : \mathbb{R} \to \mathbb{R}^2$ *bounded both in the past and in the future*

(a) $u(t) = (0,0)$ *for* $t \in \mathbb{R}$,

(b) $u(t) = (\widehat{u}_1(t), 0)$ *for* $t \in \mathbb{R}$, *where* $\widehat{u}_1$ *is the unique trajectory of* $u_1' = u_1(a_1(t) - b_{11}(t)u_1(t))$ *separated from zero and infinity given by Lemma 5.1,*

(c) $u(t) = (u_1(t), 0)$ *where* $\lim_{t \to -\infty} u_1(t) = 0$ *and* $\lim_{t \to \infty}(u_1(t) - \widehat{u}_1(t)) = 0$. *Any solution with the initial data* $(u_1(t_0), 0)$ *satisfying* $0 < u_1(t_0) < \widehat{u}_1(t_0)$ *constitutes such trajectory.*

*All solutions other than the ones named in (a)–(c) are backwards unbounded. Moreover,*

(d) *Any solution* $u(t) = (u_1(t), u_2(t))$ *with initial data* $u_1(t_0) > 0$ *and* $u_2(t_0) > 0$ *satisfies* $\lim_{t \to \infty}(u_1(t) - \widehat{u}_1(t), u_2(t)) = (0,0)$ *and* $\lim_{t \to -\infty}|u(t)| = \infty$.

(e) *Any solution* $u(t) = (0, u_2(t))$ *with* $u_2(t_0) > 0$ *satisfies* $\lim_{t \to \infty} u_2(t) = 0$ *and* $\lim_{t \to -\infty} u_2(t) = \infty$,

(f) *Any solution* $u(t) = (u_1(t), 0)$ *with* $u_1(t_0) > \widehat{u}_1(t_0)$ *satisfies* $\lim_{t \to \infty}(u_1(t) - \widehat{u}_1(t)) = 0$ *and* $\lim_{t \to -\infty} u_1(t) = \infty$.

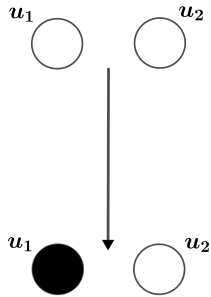

Figure 2: Two complete solutions and their connections for the two dimensional problem for which the state with existence of one species is globally asymptotically stable.

### 6.3 Extinction of both species

The last situation is the extinction of both species. To obtain this case we make the following assumptions

$$
\begin{cases}
(c_1 b_{11} + c_2 b_{12})^L > 0, \\
(c_2 b_{22} + c_1 b_{21})^L > 0,
\end{cases}
\quad
\begin{cases}
a_1 \leqslant b_{12}(d_2 + \theta c_2) - \varepsilon, \\
a_2 \leqslant b_{21}(d_1 + \theta c_1) - \varepsilon,
\end{cases}
\tag{20}
$$

with some constants $d_1, d_2, c_1, c_2 > 0$ and $\varepsilon, \theta > 0$. We observe that the first inequalities is condition $(H1)$ and the second ones is condition $(B)$ with $J = \{1, 2\}$ and $I = \varnothing$.

**Theorem 6.8.** *Assuming* (20) *the only trajectory of system* (LV-2) *which is bounded both in the past and and the future is* $u(t) = (0, 0)$ *for* $t \in \mathbb{R}$. *Moreover for every solution with initial data* $u_1(t_0) > 0$ *or* $u_2(t_0) > 0$ *we have* $\lim_{t \to \infty}(u_1(t), u_2(t)) = (0, 0)$.

*Proof.* The convergence $\lim_{t \to \infty}(u_1(t), u_2(t)) = (0, 0)$ follows from Lemma 4.2. It remains to prove that any solution with initial data $u_1(t_0) > 0$ or $u_2(t_0) > 0$ is backward unbounded, done it in Lemma 4.3. □

## 7 Attractor for non-autonomous Lotka-Volterra 3-D system

In this section we will study different possibilities of the non-autonomous attractor structure for the following three-dimensional system.

$$
\begin{cases}
u_1' = u_1(a_1(t) - b_{11}(t)u_1 - b_{12}(t)u_2 - b_{13}(t)u_3), \\
u_2' = u_2(a_2(t) - b_{21}(t)u_1 - b_{22}(t)u_2 - b_{23}(t)u_3), \\
u_3' = u_3(a_3(t) - b_{31}(t)u_1 - b_{32}(t)u_2 - b_{33}(t)u_3).
\end{cases}
\tag{LV-3}
$$

with $b_{ii}^L > 0$ and $b_{ij} \leqslant 0$ for $i, j \in \{1, 2, 3\}$. Remind that $S(t, \tau)$ for $t \geqslant \tau$ denotes the mapping that assigns to the initial condition taken at time $\tau$ the value of the solution at time $t$, now associated to (LV-3).

### 7.1 Structure of attractor for the case of permanence

The first case, that of permanence, corresponds to the dynamics depicted in Fig. 3. Following Theorem 2.1 we have to impose the following conditions on the coefficients of the problem to obtain the permanence of the three species, as well as the admissibility of all intermediate transient states together with the existence of all appropriate connections.

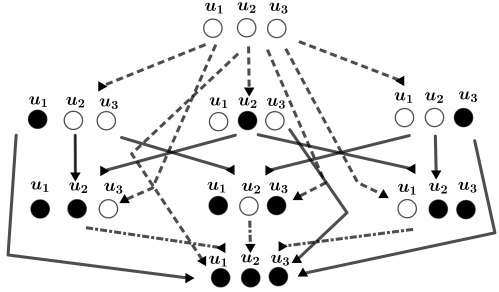

Figure 3: Eight complete solutions and their connections for the three dimensional problem with globally asymptotically stable state with coexistence for all the species.

To this end, we assume that there exist constants $c_i > 0$, $d_i > 0$ and $\bar{d}_i > 0$ for $i \in \{1, 2, 3\}$, such that the following two systems of inequalities are satisfied.

$$\left\{ \begin{array}{l} (\bar{c}_1 b_{11} + \bar{c}_2 b_{21} + \bar{c}_3 b_{31})^L > 0, \\ (\bar{c}_2 b_{22} + \bar{c}_1 b_{12} + \bar{c}_3 b_{32})^L > 0, \\ (\bar{c}_3 b_{33} + \bar{c}_1 b_{13} + \bar{c}_2 b_{23})^L > 0. \end{array} \right. \qquad \left\{ \begin{array}{l} b_{11}\bar{d}_1 \leqslant a_1 \leqslant b_{11}d_1 + b_{12}d_2 + b_{13}d_3 \\ b_{22}\bar{d}_2 \leqslant a_2 \leqslant b_{21}d_1 + b_{22}d_2 + b_{23}d_3 \\ b_{33}\bar{d}_3 \leqslant a_3 \leqslant b_{31}d_1 + b_{32}d_2 + b_{33}d_3 \end{array} \right. \qquad (21)$$

The above inequalities are exactly conditions $(A)$ and $(H2)$, respectively, for the three-dimensional case. Moreover we need to assume that

$$|b_{ij}|^U < \infty \text{ and } a_i^U < \infty \text{ for every } i, j \in \{1, 2, 3\}. \qquad (22)$$

As a direct consequence of Theorem 2.1, analogously as in the two-dimensional case, we have the following lemma.

**Lemma 7.1.** *There exists a function $u^* = (u_1^*, u_2^*, u_3^*)$ defined for $t \in \mathbb{R}$, the complete solution of* (LV-3) *such that $\bar{d}_i \leqslant u_i^* \leqslant d_i$ for $i = 1, 2, 3$. This is a unique complete trajectory bounded away from zero and infinity in all the variables. Moreover, there exist functions $(\hat{u}_1, \hat{u}_2, 0)$, $(0, \hat{u}_2, \hat{u}_3)$, $(\hat{u}_1, 0, \hat{u}_3)$, $(\bar{u}_1, 0, 0)$, $(0, \bar{u}_2, 0)$, $(0, 0, \bar{u}_3)$, and $u(t) = (0, 0, 0)$ defined for $t \in \mathbb{R}$ that are complete solutions to* (LV-3), *such that $\bar{d}_i \leqslant \hat{u}_i, \bar{u}_i \leqslant d_i$ for $i = 1, 2, 3$.*

To obtain the attractor structure we prove the existence of the connections between the trajectories obtained in the above theorem, which play the role of non-autonomous equilibria. Next theorem is our first result in this context.

**Theorem 7.2.** *There exists the trajectory of* (LV-3) *denoted by $z = (z_1, z_2, z_3) : \mathbb{R} \to \mathbb{R}^3$ such that*

$$\lim_{s \to -\infty} |(z_1(s), z_2(s), z_3(s)) - (\hat{u}_1(s), \hat{u}_2(s), 0)| = 0.$$

*and*

$$\lim_{s \to \infty} |(z_1(s), z_2(s), z_3(s)) - (u_1^*(s), u_2^*(s), u_3^*(s))| = 0.$$

*Analogous result holds for $(0, \widehat{u}_2, \widehat{u}_3)$ and $(\widehat{u}_1, 0, \widehat{u}_3)$. Furthermore, there exist a trajectory denoted by $w = (w_1, w_2, w_3) : \mathbb{R} \to \mathbb{R}^3$ such that*

$$\lim_{s \to -\infty} |(w_1(s), w_2(s), w_3(s)) - (\bar{u}_1(s), 0, 0)| = 0.$$

*and*

$$\lim_{s \to \infty} |(w_1(s), w_2(s), w_3(s)) - (u_1^*(s), u_2^*(s), u_3^*(s))| = 0.$$

*Analogous result holds for $(0, \bar{u}_2, 0)$ and $(0, 0, \bar{u}_3)$.*

*Proof.* The proof is analogous as for the two-dimensional case. First, we need to prove that the solution $(\widehat{u}_1, \widehat{u}_2, 0)$ is locally unstable, i.e. its non-autonomous unstable manifold

$$W^u((\widehat{u}_1, \widehat{u}_2, 0)) = \{(t, (w_1, w_2, w_3)) \ : \ \text{there exists a solution } z : \mathbb{R} \to \mathbb{R}^3 \text{ such that}$$
$$z(t) = (w_1, w_2, w_3) \ \text{ and } \ \lim_{s \to -\infty} |z(s) - (\widehat{u}_1(s), \widehat{u}_2(s), 0)| = 0\},$$

is nonempty and intersects the interior of the positive quadrant. We linearize the system around $(\widehat{u}_1(s), \widehat{u}_2(s), 0)$. So we write

$$w(t) = u(t) - \widehat{u}(t) \tag{23}$$

where $\widehat{u}(t) = (\widehat{u}_1(t), \widehat{u}_2(t), 0)$, and now the linearized system has the form

$$v'(t) = \begin{pmatrix} A(t) & C(t) \\ 0 & B(t) \end{pmatrix} v(t). \tag{24}$$

where

$$A(t) = \begin{pmatrix} a_1(t) - 2b_{11}(t)\widehat{u}_1(t) - b_{12}(t)\widehat{u}_2(t) & -b_{12}(t)\widehat{u}_1(t) \\ -b_{21}(t)\widehat{u}_2(t) & a_2(t) - b_{21}(t)\widehat{u}_1(t) - 2b_{22}(t)\widehat{u}_2(t) \end{pmatrix},$$

$$C(t) = \begin{pmatrix} -b_{13}(t)\widehat{u}_1(t) \\ -b_{23}(t)\widehat{u}_2(t) \end{pmatrix}$$

and $B(t) = a_3(t) - b_{31}(t)\widehat{u}_1(t) - b_{32}(t)\widehat{u}_2(t)$. We observe that the one dimensional system $v_3' = B(t)v_3$ has an exponential dichotomy with projection $P(t) = 0$ since by (21) it holds $B(t) \geqslant \delta > 0$.

Now $(\widehat{u}_1(t), \widehat{u}_2(t))$ is the solution of the two-dimensional system obtained by taking $u_3 \equiv 0$, and separated from zero and infinity. Hence, by Lemma 3.3, the system $(v_1, v_2)' = A(t)(v_1, v_2)$ has an exponential dichotomy with $P(t) = I_{2 \times 2}$.

Analogously as in two-dimensional case, we follow the results of Section 3, and since by (22) and Lemma 7.1 the function $C(t)$ is bounded, and the system (24) admits an exponential dichotomy. So by (7), the projections are given by

$$P^+ = \begin{pmatrix} 1 & 0 & 0 \\ 0 & 1 & 0 \\ 0 & 0 & 0 \end{pmatrix} \text{ on } \mathbb{R}^+ \ \text{ and } P^- = \begin{pmatrix} 1 & 0 & L_1^- \\ 0 & 1 & L_2^- \\ 0 & 0 & 0 \end{pmatrix} \text{ on } \mathbb{R}^-.$$

Now the projection

$$P = P^- = \begin{pmatrix} 1 & 0 & L_1^- \\ 0 & 1 & L_2^- \\ 0 & 0 & 0 \end{pmatrix},$$

has the same range as $P^+$ and hence the system (24) has an exponential dichotomy with $P$ on whole $\mathbb{R}$, and

$$P(t) = \begin{pmatrix} 1 & 0 & L_1(t) \\ 0 & 1 & L_2(t) \\ 0 & 0 & 0 \end{pmatrix},$$

for certain bounded functions $L_1(t)$ and $L_2(t)$. Note that range$(P(t)) = \{(\beta_1, \beta_2, 0) : \beta_1, \beta_2 \in \mathbb{R}\}$, of dimension 2, is always the stable space in the linearization, while

$$\text{range}(I - P(t)) = \left\{ \alpha \begin{pmatrix} -L_1(t) \\ -L_2(t) \\ 1 \end{pmatrix} : \alpha \in \mathbb{R} \right\},$$

of dimension 1, is the time dependent unstable space in the linearization.

Analogously as in the two-dimensional case we use the non-autonomous unstable manifold theorem, cf., [KR11, Theorem 6.10 and Exercise 6.11]. There exists a neighborhood $U$ of zero in $\mathbb{R}^3$ and a function $\Sigma^- : \mathbb{R} \times U \to \mathbb{R}^3$ such that $\Sigma^-(t, x) = \Sigma^-(t, (I - P(t))x) \in \text{range}(P(t))$, $\Sigma^-(t, 0) = 0$, $\lim_{|x| \to 0} \frac{\Sigma^-(t,x)}{|x|} = 0$. Moreover, if for some $t \in \mathbb{R}$ and $x \in \text{range}(I - P(t)) \cap U$ we have $y = x + \Sigma^-(t, x)$, then, provided $|y|$ is small enough, the backward trajectory $w : (-\infty, t] \to \mathbb{R}^3$ of (23) with $w(t) = y$ belongs to the graph of $\Sigma^-$, i.e. $\Sigma^-(\tau, (I - P(\tau))w(\tau)) = P(\tau)w(\tau)$ and satisfies $|w(\tau)| \leq C e^{\delta(\tau - t)}$ for every $\tau \leq t$. Moreover for a small positive $\alpha$ the third coordinate of the point

$$\alpha \begin{pmatrix} -L_1(t) \\ -L_2(t) \\ 1 \end{pmatrix} + \Sigma^- \left( t, \alpha \begin{pmatrix} -L_1(t) \\ -L_2(t) \\ 1 \end{pmatrix} \right)$$

is equal to $\alpha$ and hence positive, while the first two coordinates are small enough so that in the original coordinates they are also positive. Hence, for any $t \in \mathbb{R}$ there exists the point with all the coordinates positive at time $t$ which is backward exponentially attracted to zero. In terms of the original system (LV-3) in view of [Red96, Theorem 1 (iii) and (vi)] this gives us the assertion of the Theorem.

The proof of local instability of $(\bar{u}_1, 0, 0)$ is analogous. The system linearized around $(\bar{u}_1, 0, 0)$ has the form

$$v'(t) = \begin{pmatrix} A(t) & C(t) \\ 0 & B(t) \end{pmatrix} v(t). \tag{25}$$

where

$$A(t) = \begin{pmatrix} a_1(t) - 2b_{11}(t)\bar{u}_1(t) \end{pmatrix}, \quad B(t) = \begin{pmatrix} a_2(t) - b_{21}(t)\bar{u}_1(t) & 0 \\ 0 & a_3(t) - b_{31}(t)\bar{u}_1(t) \end{pmatrix}$$

$$C(t) = \begin{pmatrix} -b_{12}(t)\bar{u}_1(t) & -b_{13}(t)\bar{u}_1(t) \end{pmatrix}$$

Now, the two dimensional system $(v_2(t), v_3(t))' = B(t)(v_2(t), v_3(t))$ has an exponential dichotomy with projection $P(t) = 0$ since by (21) both diagonal terms are positive and separated from zero. The system $v_1' = A(t)v_1$, by Theorem 6.3, has an exponential dichotomy with projection $P(t) = I$, We follows the same process as the first case: the function $C(t)$ is bounded, and the system (24) admits an exponential dichotomy, such that by (7) the projections are given by

$$P^+ = \begin{pmatrix} 1 & 0 & 0 \\ 0 & 0 & 0 \\ 0 & 0 & 0 \end{pmatrix} \text{ on } \mathbb{R}^+ \text{ and } P^- = \begin{pmatrix} 1 & L_1^- & L_2^- \\ 0 & 0 & 0 \\ 0 & 0 & 0 \end{pmatrix} \text{ on } \mathbb{R}^-.$$

Since the projection

$$P = P^- = \begin{pmatrix} 1 & L_1^- & L_2^- \\ 0 & 0 & 0 \\ 0 & 0 & 0 \end{pmatrix},$$

has the same range as $P^+$, hence the system (24) has exponential dichotomy on $\mathbb{R}$ with the projection

$$P(t) = \begin{pmatrix} 1 & L_1(t) & L_2(t) \\ 0 & 0 & 0 \\ 0 & 0 & 0 \end{pmatrix},$$

for certain bounded functions $L_1(t)$ and $L_2(t)$. Note that in this case range$(P(t)) = \{\beta(1, 0, 0) : \beta \in \mathbb{R}\}$, is of dimension 1, and is always the stable space in the linearization, while

$$\text{range}(I - P(t)) = \left\{ \alpha \begin{pmatrix} -L_1(t) \\ 1 \\ 0 \end{pmatrix} + \beta \begin{pmatrix} -L_2(t) \\ 0 \\ 1 \end{pmatrix} : \alpha, \beta \in \mathbb{R} \right\},$$

is a two dimensional space,and it is the time dependent unstable space in the linearization. The end of the proof is analogous to the first case. with these new spaces. $\qquad\square$

**Theorem 7.3.** *There exists a trajectory of* (LV-3) *denoted by* $y = (y_1, y_2, y_3) : \mathbb{R} \to \mathbb{R}^3$ *such that*

$$\lim_{s \to -\infty} |(y_1(s), y_2(s), y_3(s)| = 0.$$

*and*

$$\lim_{s \to \infty} |(y_1(s), y_2(s), y_3(s)) - (u_1^*(s), u_2^*(s), u_3^*(s))| = 0.$$

*Proof.* The proof is analogous the the proof in the two-dimensional case, cf. Theorem 6.4, with following linearized system at $(0, 0, 0)$

$$\begin{cases} v'(t) = \begin{pmatrix} a_1(t) & 0 & 0 \\ 0 & a_2(t) & 0 \\ 0 & 0 & a_3(t) \end{pmatrix} v(t), \end{cases} \tag{26}$$

$\qquad\square$

We prove now the result of the three-dimensional analogy of Theorem 6.5. The proof follows the same argument as in Theorem 6.5.

**Theorem 7.4.** *If a trajectory with a strictly positive initial condition is bounded both in the past and in the future then it must be one of trajectories described above.*

*Proof.* From [Red96, Theorem 1. (v)] we know that $u^*$ is the only solution on $\mathbb{R}$ that is bounded away from 0 and $\infty$. Thus, if the complete trajectory $u = (u_1, u_2, u_3)$ is bounded, not coincides with $u^*$, and none of the three coordinates is zero, then for at least one of its coordinates there must exist a decreasing sequence $t_n \to -\infty$ such that $u_i(t_n) \to 0$ and $u_i(t_n) < \overline{d}_i$. We consider two cases.

*Case 1. Two coordinates converge to zero backward and one is separated from zero.* Assume that there exist two sequences $t_n \to -\infty$ and $\hat{t}_n \to -\infty$ such that $u_1(t_n) \to 0$, $u_2(\hat{t}_n) \to 0$, $u_1(t_n) < \overline{d}_1$ and $u_2(\hat{t}_n) < \overline{d}_2$ We first prove that $\lim_{t \to -\infty} u_i(t) = 0$ for $i = 1, 2$. We proceed analogously as in Lemma 5.2 and Theorem 6.5. If for some $t \in (t_{n+1}, t_n)$, we have $u_1(t) > u_1(t_n)$ then there exists $t^* \in [t, t_n)$ such that $u_1(t^*) > u_1(t_n)$ and on $(t^*, t_n)$ the function $u_1$ is strictly less that $\overline{d}_1$. Hence

$$0 > u_1(t_n) - u_1(t^*) = \int_{t^*}^{t_n} u_1(s)(a_1(s) - b_{11}(s)u_1(s) - b_{12}(s)u_2(s) - b_{13}u_3(s))\, ds \geqslant$$

$$\geqslant \int_{t^*}^{t_n} u_1(s)(a_1(s) - b_{11}(s)\overline{d}_1)ds \geqslant 0,$$

and we have the contradiction. The same argument allows us to get $\lim_{t \to -\infty} u_2(t) = 0$.

Now we need to prove that for $u_3(t)$, that is separated from zero, it must be

$$\lim_{t \to -\infty} |u_3(t) - \bar{u}_3(t)| = 0.$$

To this end assume that for some sequence $t_n \to -\infty$

$$\lim_{n \to \infty} |u_3(t_n) - \bar{u}_3(t_n)| = c > 0.$$

Now, since $u_3(t)$ is separated from zero, denoting $K = (u_3)^L$ and $M = (u_3)^U$ and using Lemma 2.2 and Lemma 2.3, we obtain for every $x \in [K, M]$

$$|u_3(t_n) - \bar{u}_3(t_n)| \leqslant |S(t_n, t_n - t)(u_1(t_n - t), u_2(t_n - t), u_3(t_n - t)) - S(t_n, t_n - t)(0, 0, x)|$$
$$+ |S(t_n, t_n - t)(0, 0, x) - S(t_n, t_n - t)(0, 0, \bar{u}_3(t_n - t))|$$
$$\leqslant e^{\kappa t}|(u_1(t_n - t), u_2(t_n - t), u_3(t_n - t)) - (0, 0, x)| + \frac{2 \max\{M, d_3\}^2}{\min\{K, \overline{d}_3\}} e^{-\delta \min\{K, \overline{d}_3\}t},$$

with a constant $\kappa$ depending on the bounds on $u$ and the coefficients of the problem. Pick $\varepsilon > 0$. We can find $t > 0$ such that $\frac{2 \max\{M, d_3\}^2}{\min\{K, \overline{d}_3\}} e^{-\delta \min\{K, \overline{d}_3\}t} < \frac{\varepsilon}{2}$. The sequence $u_3(t_n - t)$ has a convergent subsequence, let $x$ be a limit of this subsequence and pick $n$ large enough such that $e^{\kappa t}|(u_1(t_n - t), u_2(t_n - t), u_3(t_n - t)) - (0, 0, x)| < \frac{\varepsilon}{2}$. Hence for $n$ large enough, on a subsequence,

$$|u_3(t_n) - \bar{u}_3(t_n)| < \varepsilon,$$

a contradiction.

*Case 2. Only one coordinate converges to zero backwards.* We assume, for example, that there exists a sequence $t_n \to -\infty$ such that $u_1(t_n) \to 0$, then, analogously to the previous case, we can prove that $\lim_{t \to -\infty} u_1(t) = 0$. We must prove that for $u_2(t)$, and $u_3(t)$, that are backward separated from zero, it must be

$$\lim_{t \to -\infty} |(u_2(t), u_3(t)) - (\hat{u}_2(t), \hat{u}_3(t))| = 0.$$

To this end assume that for some sequence $t_n \to -\infty$

$$\lim_{n \to \infty} |(u_2(t_n), u_3(t_n)) - (\widehat{u}_2(t_n), \widehat{u}_3(t_n))| = c > 0. \tag{27}$$

For every $x_i \in [u_i^L, u_i^U]$ for $i \in \{2, 3\}$, and using Lemma 2.3 with $(H2)$ and Lemma 2.2 we have for $t \geqslant 0$

$$\begin{aligned}
|(u_2(t_n), &u_3(t_n)) - (\widehat{u}_2(t_n), \widehat{u}_3(t_n))| \\
&\leqslant |S(t_n, t_n - t)(u_1(t_n - t), u_2(t_n - t), u_3(t_n - t)) - S(t_n, t_n - t)(0, x_2, x_3)| \\
&\quad + |S(t_n, t_n - t)(0, x_2, x_3) - S(t_n, t_n - t)(0, \widehat{u}_2(t_n - t), \widehat{u}_3(t_n - t))| \\
&\leqslant e^{\kappa t} |(u_1(t_n - t), u_2(t_n - t), u_3(t_n - t)) - (0, x_2, x_3)| + \frac{2\sigma^2}{\bar{\sigma}} e^{-\bar{\sigma}\delta t},
\end{aligned}$$

where $\kappa, \sigma, \bar{\sigma}$ depend on the bounds on coefficients of the problem, constants appearing in (21) and on $u_i^L, u_i^U$. In particular [Red96, Lemma 2] implies that the solution starting from $(0, x_2, x_3)$ at any time is separated from zero by a constant independent on initial time and the choice of initial data as long as $x_i \in [u_i^L, u_i^U]$. Pick $\varepsilon > 0$. We can find $t > 0$ such that $\frac{2\sigma^2}{\bar{\sigma}} e^{-\bar{\sigma}\delta t} < \frac{\varepsilon}{2}$. The sequence $(u_2(t_n - t), u_3(t_n - t))$ has a convergent subsequence, let $(x_2, x_3)$ be a limit of this subsequence and associated times are denoted by $\{\widehat{t}_n\}$. We pick $n$ large enough such that $e^{\kappa t} |(u_1(t_n - t), u_2(t_n - t), u_3(t_n - t)) - (0, x_2, x_3)| < \frac{\varepsilon}{2}$, which leads into a direct contradiction with (27). □

Finally, the full image of the attractor represented in Figure 3 is described in the following theorem with a summary of the results obtained in this section and using Lemma 5.2, Theorem 6.3 and Theorem 6.4 for the remaining heterolinics connections.

**Theorem 7.5.** *Assuming* (21) *the system* (LV-3) *has the following solutions* $u : \mathbb{R} \to \mathbb{R}^3$ *bounded both in the past and in the future:*

(a) $u(t) = (0, 0, 0)$ *for* $t \in \mathbb{R}$,

(b) $u(t) = (\bar{u}_1(t), 0, 0), u(t) = (0, \bar{u}_2(t), 0)$ *and* $u(t) = (0, 0, \bar{u}_3)$, *given in Lemma 7.1.*

(c) *Any solution* $u(t) = (u_1(t), 0, 0)$ *with initial condition* $0 < u_1(t_0) < \bar{u}_1(t_0)$, *where* $\lim_{t \to -\infty} u_1(t) = 0$ *and* $\lim_{t \to \infty} (u_1(t) - \bar{u}_1(t)) = 0$. *Analogously,* $u(t) = (0, u_2(t), 0)$ *with initial condition* $0 < u_2(t_0) < \bar{u}_2(t_0)$, *where* $\lim_{t \to -\infty} u_2(t) = 0$ *and* $\lim_{t \to \infty} (u_2(t) - \bar{u}_2(t)) = 0$, *and* $u(t) = (0, 0, u_3(t))$ *with initial condition* $0 < u_3(t_0) < \bar{u}_3(t_0)$, *where* $\lim_{t \to -\infty} u_3(t) = 0$ *and* $\lim_{t \to \infty} (u_3(t) - \bar{u}_3(t)) = 0$.

(d) $u(t) = (\widehat{u}_1(t), \widehat{u}_2(t), 0), \ u(t) = (\widehat{u}_1(t), 0, \widehat{u}_3(t))$ *and* $u(t) = (0, \widehat{u}_2(t), \widehat{u}_3(t))$, *given in Lemma 7.1.*

(e) $u(t) = (u_1(t), u_2(t), 0)$ *where* $\lim_{t \to -\infty} (u_1(t), u_2(t)) = (0, 0)$ *and* $\lim_{t \to \infty} ((u_1(t), u_2(t)) - (\widehat{u}_1(t), \widehat{u}_2(t))) = 0$. *Analogously,* $u(t) = (u_1(t), 0, u_3(t))$ *where* $\lim_{t \to -\infty} (u_1(t), u_3(t)) = (0, 0)$ *and* $\lim_{t \to \infty} ((u_1(t), u_3(t)) - (\widehat{u}_1(t), \widehat{u}_3(t))) = 0$, *and* $u(t) = (0, u_2(t), u_3(t))$ *where* $\lim_{t \to -\infty} (u_2(t), u_3(t)) = (0, 0)$ *and* $\lim_{t \to \infty} ((u_2(t), u_3(t)) - (\widehat{u}_2(t), \widehat{u}_3(t))) = 0$.

(f) $u(t) = (u_1(t), u_2(t), 0)$ where $\lim_{t\to-\infty}(u_1(t), u_2(t)) = (\bar{u}_1, 0)$ and $\lim_{t\to\infty}((u_1(t), u_2(t)) - (\widehat{u}_1(t), \widehat{u}_2(t))) = 0$. Analogously, $u(t) = (u_1(t), u_2(t), 0)$ where $\lim_{t\to-\infty}(u_1(t), u_2(t)) = (0, \bar{u}_2)$ and $\lim_{t\to\infty}((u_1(t), u_2(t)) - (\widehat{u}_1(t), \widehat{u}_2(t))) = 0$.

(g) $u(t) = (u_1(t), 0, u_3(t))$ where $\lim_{t\to-\infty}(u_1(t), u_3(t)) = (\bar{u}_1, 0)$ and $\lim_{t\to\infty}((u_1(t), u_3(t)) - (\widehat{u}_1(t), \widehat{u}_3(t))) = 0$. Analogously, $u(t) = (u_1(t), 0, u_3(t))$ where $\lim_{t\to-\infty}(u_1(t), u_3(t)) = (0, \bar{u}_3)$ and $\lim_{t\to\infty}((u_1(t), u_3(t)) - (\widehat{u}_1(t), \widehat{u}_3(t))) = 0$.

(h) $u(t) = (0, u_2(t), u_3(t))$ where $\lim_{t\to-\infty}(u_2(t), u_3(t)) = (\bar{u}_2, 0)$ and $\lim_{t\to\infty}((u_2(t), u_3(t)) - (\widehat{u}_2(t), \widehat{u}_3(t))) = 0$. Analogously, $u(t) = (0, u_2(t), u_3(t))$ where $\lim_{t\to-\infty}(u_2(t), u_3(t)) = (0, \bar{u}_3)$ and $\lim_{t\to\infty}((u_2(t), u_3(t)) - (\widehat{u}_2(t), \widehat{u}_3(t))) = 0$.

(i) $u(t) = (u_1^*(t), u_2^*(t), u_3^*(t))$ bounded away from zero and infinity given in Lemma 7.1.

(j) $u(t) = (u_1(t), u_2(t), u_3(t))$ such that $\lim_{t\to-\infty}(u_1(t), u_2(t), u_3(t)) = (0, 0, 0)$ and $\lim_{t\to\infty}(u_1(t) - u_1^*(t), u_2(t) - u_2^*(t), u_3(t) - u_3^*(t)) = 0$, given in Theorem 7.3.

(k) $u(t) = (u_1(t), u_2(t), u_3(t))$ such that $\lim_{t\to-\infty}(u_1(t), u_2(t), u_3(t)) = (\widehat{u}_1(t), \widehat{u}_2(t), 0)$ and $\lim_{t\to\infty}(u_1(t) - u_1^*(t), u_2(t) - u_2^*(t), u_3(t) - u_3^*(t)) = 0$, given in Theorem 7.2. Analogously, $u(t) = (u_1(t), u_2(t), u_3(t))$ such that $\lim_{t\to-\infty}(u_1(t), u_2(t), u_3(t)) = (\widehat{u}_1(t), 0, \widehat{u}_3(t))$ and $\lim_{t\to\infty}(u_1(t) - u_1^*(t), u_2(t) - u_2^*(t), u_3(t) - u_3^*(t)) = 0$, and $u(t) = (u_1(t), u_2(t), u_3(t))$ such that $\lim_{t\to-\infty}(u_1(t), u_2(t), u_3(t)) = (0, \widehat{u}_2(t), \widehat{u}_3(t))$ and $\lim_{t\to\infty}(u_1(t) - u_1^*(t), u_2(t) - u_2^*(t), u_3(t) - u_3^*(t)) = 0$.

(l) $u(t) = (u_1(t), u_2(t), u_3(t))$ such that $\lim_{t\to-\infty}(u_1(t), u_2(t), u_3(t)) = (\bar{u}_1(t), 0, 0)$ and $\lim_{t\to\infty}(u_1(t) - u_1^*(t), u_2(t) - u_2^*(t), u_3(t) - u_3^*(t))$, given in Theorem 7.2. Analogously, $u(t) = (u_1(t), u_2(t), u_3(t))$ such that $\lim_{t\to-\infty}(u_1(t), u_2(t), u_3(t)) = (0, \bar{u}_2(t), 0)$ and $\lim_{t\to\infty}(u_1(t) - u_1^*(t), u_2(t) - u_2^*(t), u_3(t) - u_3^*(t)) = 0$, and $u(t) = (u_1(t), u_2(t), u_3(t))$ such that $\lim_{t\to-\infty}(u_1(t), u_2(t), u_3(t)) = (0, 0, \bar{u}_3(t))$ and $\lim_{t\to\infty}(u_1(t) - u_1^*(t), u_2(t) - u_2^*(t), u_3(t) - u_3^*(t)) = 0$.

*Moreover, any solution of* (LV-3) *bounded both in the past and in the future is one of the solutions described in items (a)-(l).*

## 7.2  Structure of attractor for the case of extinction of one species

In case of the permanence of the two species we would have a scheme depicted in the next Figure.

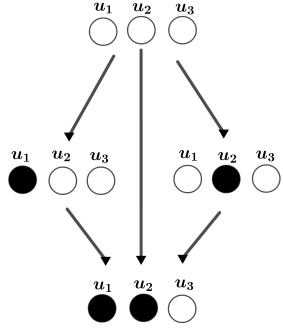

Figure 4: Four complete solutions and their connections for the three dimensional problem with globally asymptotically stable state with coexistence of two of the three species.

Following Theorem 4.5 we have to impose the following conditions to $a_1, a_2$ and $a_3$ to obtain the extinction of one species.

$$\left\{ \begin{array}{l} (c_1 b_{11} + c_2 b_{12} + c_3 b_{13})^L > 0, \\ (c_2 b_{22} + c_1 b_{21} + c_3 b_{23})^L > 0, \\ (c_3 b_{33} + c_1 b_{31} + c_2 b_{32})^L > 0, \end{array} \right. \quad \left\{ \begin{array}{l} b_{11}\bar{d}_1 + \varepsilon \leqslant a_1 \leqslant b_{11}d_1 + b_{12}d_2 + b_{13}d_3 + \theta(b_{12}c_2 + b_{13}c_3) - \varepsilon, \\ b_{22}\bar{d}_2 + \varepsilon \leqslant a_2 \leqslant b_{21}d_1 + b_{22}d_2 + b_{23}d_3 + \theta(b_{21}c_1 + b_{23}c_3) - \varepsilon, \\ a_3 \leqslant b_{31}d_1 + b_{32}d_2 + \theta(b_{31}c_1 + b_{32}c_2) - \varepsilon, \end{array} \right. \tag{28}$$

with some constants $d_i, c_i > 0$ for $i = 1, 2, 3$ , $\bar{d}_j > 0$ for $j = 1, 2$ and $\varepsilon, \theta > 0$. The above inequalities are the conditions $(H1)$ and $(B2)$ for the three-dimensional case. Moreover, we need the condition $(H2)$ for the two dimensional subproblem, namely

$$\left\{ \begin{array}{l} (\bar{c}_1 b_{11} + \bar{c}_2 b_{21})^L > 0, \\ (\bar{c}_2 b_{22} + \bar{c}_1 b_{12})^L > 0, \end{array} \right. \tag{29}$$

to obtain the existence of a solution $(\hat{u}_1, \hat{u}_2)$ bounded away from zero and infinity attracting given by Lemma 6.1. Also we need to assume that $|b_{ij}|^U < \infty$ and $|a_i|^U < \infty$ for every $i \in \{1, 2\}$ and $j \in \{1, 2, 3\}$. Summarizing all previous results we can formulate the following theorem.

**Theorem 7.6.** *Assuming* (28) *and* (29) *the system* (LV-3) *has the following trajectories* $u : \mathbb{R} \to \mathbb{R}^3$ *bounded both in the past and in the future*

(a) $u(t) = (0, 0, 0)$ *for* $t \in \mathbb{R}$,

(b) $u(t) = (\bar{u}_1(t), 0, 0), u(t) = (0, \bar{u}_2(t), 0)$ *for* $t \in \mathbb{R}$, *where* $\bar{u}_i$, *for* $i = 1, 2$, *is the unique trajectory of* $u'_i = u_i(a_i(t) - b_{ii}(t)u_i(t))$ *separated from zero and infinity given by Lemma 5.1,*

(c) $u(t) = (u_1(t), 0, 0)$ where $\lim_{t \to -\infty} u_1(t) = 0$ and $\lim_{t \to \infty}(u_1(t) - \bar{u}_1(t)) = 0$. Any solution with the initial data $(u_1(t_0), 0, 0)$ satisfying $0 < u_1(t_0) < \bar{u}_1(t_0)$ constitutes such trajectory. Analogously with any solution with the initial data $(0, u_2(t_0), 0)$ satisfying $0 < u_2(t_0) < \bar{u}_2(t_0)$, constitutes the trajectory that satisfies $u(t) = (0, u_2(t), 0)$ where $\lim_{t \to -\infty} u_2(t) = 0$ and $\lim_{t \to \infty}(u_2(t) - \bar{u}_2(t)) = 0$.

(d) $u(t) = (\hat{u}_1(t), \hat{u}_2(t), 0)$, where $(\hat{u}_1, \hat{u}_2)$, is the unique trajectory of the two dimensional system obatined taking $u_3 \equiv 0$ separated from zero and infinity given by Lemma 6.1.

(e) $u(t) = (u_1(t), u_2(t), 0)$ where $\lim_{t \to -\infty}(u_1(t), u_2(t)) = (0, 0)$ and $\lim_{t \to \infty}((u_1(t), u_2(t)) - (\hat{u}_1(t), \hat{u}_2(t))) = 0$.

(f) $u(t) = (u_1(t), u_2(t), 0)$ where $\lim_{t \to -\infty}(u_1(t), u_2(t)) = (\bar{u}_1, 0)$ and $\lim_{t \to \infty}((u_1(t), u_2(t)) - (\hat{u}_1(t), \hat{u}_2(t))) = 0$. Analogously, $u(t) = (u_1(t), u_2(t), 0)$ where $\lim_{t \to -\infty}(u_1(t), u_2(t)) = (0, \bar{u}_2)$ and $\lim_{t \to \infty}((u_1(t), u_2(t)) - (\hat{u}_1(t), \hat{u}_2(t))) = 0$.

If a solution $u : \mathbb{R} \to \mathbb{R}^3$ is bounded then it must be on of the solutions of items (a)-(f). Moreover,

(g) Any solution $u(t) = (u_1(t), u_2(t), u_3(t))$ with initial data $u_i(t_0) > 0$ for $i = 1, 2, 3$, satisfies $\lim_{t \to \infty}(u_1(t) - \hat{u}_1(t), u_2(t) - \hat{u}_2(t), u_3(t)) = (0, 0, 0)$ and $\lim_{t \to -\infty} |u(t)| = \infty$.

(h) Any solution $u(t) = (0, 0, u_3(t))$ with $u_3(t_0) > 0$ satisfies $\lim_{t \to \infty} u_3(t) = 0$ and $\lim_{t \to -\infty} u_3(t) = \infty$,

(i) Any solution $u(t) = (u_1(t), 0, 0)$ with $u_1(t_0) > \bar{u}_1(t_0)$ satisfies $\lim_{t \to \infty}(u_1(t) - \bar{u}_1(t)) = 0$ and $\lim_{t \to -\infty} u_1(t) = \infty$. Analogously, any solution $u(t) = (0, u_2(t), 0)$ with $u_2(t_0) > \bar{u}_2(t_0)$ satisfies $\lim_{t \to \infty}(u_2(t) - \bar{u}_2(t)) = 0$ and $\lim_{t \to -\infty} u_2(t) = \infty$.

(j) Any solution $u(t) = (u_1(t), 0, u_3(t))$ with $u_1(t_0) > 0, u_3(t_0) > 0$ satisfies $\lim_{t \to \infty}(u_1(t) - \bar{u}_1(t), u_3(t)) = 0$ and $\lim_{t \to -\infty} |(u_1(t), u_3(t))| = \infty$. Analogously, any solution $u(t) = (0, u_2(t), u_3(t))$ with $u_2(t_0) > 0, u_3(t) > 0$ satisfies $\lim_{t \to \infty}(u_2(t) - \bar{u}_2(t), u_3(t)) = 0$ and $\lim_{t \to -\infty} |(u_2(t), u_3(t))| = \infty$.

## 7.3 The case of extinction of two species

In case of the extinction of the two species we would have a scheme as the next Figure.

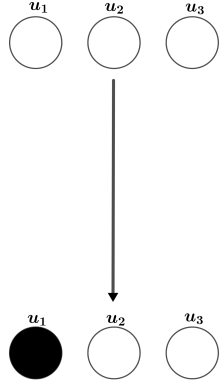

Figure 5: Two complete solutions and their connections for the three dimensional problem with globally asymptotically stable state with existence of one species.

By Theorem 4.5 we have to impose the following conditions on the problem coefficients to obtain the case of the permanence of one species and extinction of the remaining two.

$$
\left\{
\begin{array}{l}
(c_1 b_{11} + c_2 b_{12} + c_3 b_{13})^L > 0, \\
(c_2 b_{22} + c_1 b_{21} + c_3 b_{23})^L > 0, \\
(c_3 b_{33} + c_1 b_{31} + c_2 b_{32})^L > 0,
\end{array}
\right.
\qquad
\left\{
\begin{array}{l}
b_{11}\bar{d}_1 + \varepsilon \leqslant a_1 \leqslant b_{11}d_1 + b_{12}d_2 + b_{13}d_3 + \theta(b_{12}c_2 + b_{13}c_3) - \varepsilon, \\
a_2 \leqslant b_{21}d_1 + b_{23}d_3 + \theta(b_{21}c_1 + b_{23}c_3) - \varepsilon, \\
a_3 \leqslant b_{31}d_1 + b_{32}d_2 + \theta(b_{31}c_1 + b_{32}c_2) - \varepsilon,
\end{array}
\right.
$$
$$(30)$$

with some constants $d_i, c_i > 0$ for $i = 1, 2, 3$, $\bar{d}_1 > 0$ and $\varepsilon, \theta > 0$. The above inequalitites are the conditions $(H1)$ and $(B1)$ for the three-dimensional case. We need to assume again that $|b_{1j}|^U < \infty$ and $|a_1|^U < \infty$ for every $j \in \{1, 2, 3\}$

**Theorem 7.7.** *Assuming* (30) *the system* (LV-3) *has the following trajectories* $u : \mathbb{R} \to \mathbb{R}^3$ *bounded both in the past and in the future*

(a) $u(t) = (0, 0, 0)$ *for* $t \in \mathbb{R}$,

(b) $u(t) = (\bar{u}_1(t), 0, 0)$, *for* $t \in \mathbb{R}$, *where* $\bar{u}_1$ *is the unique trajectory of* $u'_1 = u_1(a_1(t) - b_{11}(t)u_1(t))$ *separated from zero and infinity given by Lemma 5.1,*

(c) $u(t) = (u_1(t), 0, 0)$ *where* $\lim_{t \to -\infty} u_1(t) = 0$ *and* $\lim_{t \to \infty}(u_1(t) - \bar{u}_1(t)) = 0$. *Any solution with the initial data* $(u_1(t_0), 0, 0)$ *satisfying* $0 < u_1(t_0) < \bar{u}_1(t_0)$ *constitutes such trajectory.*

*The solutions described in (a)-(c) are the only ones which are bounded both in the past and in the future. Moreover,*

(d) *Any solution* $u(t) = (u_1(t), u_2(t), u_3(t))$ *with initial data* $u_i(t_0) > 0$ *for* $i = 1, 2, 3$, *satisfies* $\lim_{t \to \infty}(u_1(t) - \bar{u}_1(t), u_2(t), u_3(t)) = (0, 0, 0)$ *and* $\lim_{t \to -\infty} |u(t)| = \infty$.

(e) *Any solution* $u(t) = (0, u_2(t), 0)$ *with* $u_2(t_0) > 0$ *satisfies* $\lim_{t \to \infty} u_2(t) = 0$ *and* $\lim_{t \to -\infty} u_2(t) = \infty$. *Analogously,* $u(t) = (0, 0, u_3(t))$ *with* $u_3(t_0) > 0$ *satisfies* $\lim_{t \to \infty} u_3(t) = 0$ *and* $\lim_{t \to -\infty} u_3(t) = \infty$,

(f) *Any solution* $u(t) = (u_1(t), 0, 0)$ *with* $u_1(t_0) > \bar{u}_1(t_0)$ *satisfies* $\lim_{t \to \infty}(u_1(t) - \bar{u}_1(t)) = 0$ *and* $\lim_{t \to -\infty} u_1(t) = \infty$.

(g) *Any solution* $(0, u_2(t), u_3(t))$ *with* $u_2(t_0) > 0, u_3(t_0) > 0$ *satisfies* $\lim_{t \to \infty}(u_2(t), u_3(t)) = 0$ *and* $\lim_{t \to -\infty} |(u_2(t), u_3(t))| = \infty$.

(h) *Any solution* $(u_1(t), 0, u_3(t))$ *with* $u_1(t_0) > 0$ *and* $u_3(t_0) > 0$ *satisfies* $\lim_{t \to \infty}(u_1(t) - \bar{u}_1(t), 0, u_3(t)) = 0$ *and* $\lim_{t \to -\infty} |(u_1(t), u_3(t))| = \infty$. *Analogously any solution* $(u_1(t), u_2(t), 0)$ *with* $u_1(t_0) > 0$ *and* $u_2(t_0) > 0$ *satisfies* $\lim_{t \to \infty}(u_1(t) - \bar{u}_1(t), u_2(t), 0) = 0$ *and* $\lim_{t \to -\infty} |(u_1(t), u_2(t))| = \infty$.

## 7.4 Extinction of all species

The last situation is the extinction of all species. To obtain this case we make the following assumptions

$$\begin{cases} (c_1 b_{11} + c_2 b_{12} + c_3 b_{13})^L > 0, \\ (c_2 b_{22} + c_1 b_{21} + c_3 b_{23})^L > 0, \\ (c_3 b_{33} + c_1 b_{31} + c_2 b_{32})^L > 0. \end{cases} \quad \begin{cases} a_1 \leqslant b_{12} d_2 + b_{13} d_3 + \theta(b_{12} c_2 + b_{23} c_3) - \varepsilon \\ a_2 \leqslant b_{21} d_1 + b_{23} d_3 + \theta(b_{21} c_1 + b_{23} c_3) - \varepsilon \\ a_3 \leqslant b_{31} d_1 + b_{32} d_2 + \theta(b_{31} c_1 + b_{32} c_2) - \varepsilon \end{cases} \quad (31)$$

with some constants $d_i, c_i > 0$ for $i \in \{1, 2, 3\}$, and $\varepsilon, \theta > 0$. We observe that the first inequalities is the condition $(H1)$ and the second ones is the condition $(B)$ with $J = \{1, 2, 3\}$ and $I = \varnothing$.

**Theorem 7.8.** *Assuming* (31) *the only trajectory of system* (LV-3) *which is bounded both in the past and and the future is* $u(t) = (0, 0, 0)$ *for* $t \in \mathbb{R}$. *Moreover for every solution with initial data* $u_i(t_0) > 0$ *for any* $i \in \{1, 2, 3\}$, *we have* $\lim_{t \to \infty}(u_1(t), u_2(t), u_3(t)) = (0, 0, 0)$ *and* $\lim_{t \to -\infty} |(u_1(t), u_2(t), u_3(t))| = \infty$.

*Proof.* The convergence $\lim_{t \to \infty}(u_1(t), u_2(t), u_3(t)) = (0, 0, 0)$ follows from Lemma 4.2. The fact that any solution with initial data $u_1(t_0) > 0$, $u_2(t_0)$ or $u_3(t_0) > 0$ is backward unbounded, follows from Lemma 4.3. $\square$

## Funding

Work of JGF has been partially supported by the Spanish Ministerio de Ciencia, Innovación y Universidades (MCIU), Agencia Estatal de Investigación (AEI), Fondo Europeo de Desarrollo Regional (FEDER) under grant PRE2019-087385. JGF and JALR have been also supported by projects PGC2018-096540-B-I00 and PID2021-122991NB-C21. Work of PK has been supported by Polish National Agency for Academic Exchange (NAWA) within the Bekker Programme under Project No. PPN/BEK/2020/1/00265/U/00001, and by National Science Center (NCN) of Poland under Projects No. UMO-2016/22/A/ST1/00077 and DEC-2017/25/B/ST1/00302. AS has been partially supported by projects US-1381261 and PGC2018-098308-B-I00.

# References

[Ahm93] S. Ahmad. On the nonautonomous Volterra–Lotka competition equations. *Proceedings of the American Mathematical Society*, 117:199–204, 1993.

[AKL22] P. Almaraz, P. Kalita, and J.A. Langa. Structural stability of invasion graphs for generalized lotka–volterra systems. https://arxiv.org/abs/2209.09802v4, 2022.

[AL95] S. Ahmad and A. Lazer. On the nonautonomous $n$-competing species problems. *Applicable Analysis.*, 57:309–323, 1995.

[AL98] S. Ahmad and A. Lazer. Necessary and sufficient average growth in a Lotka–Volterra system. *Nonlinear Analysis*, 34:191–228, 1998.

[AL00] S. Ahmad and A. Lazer. Average conditions for global asymptotic stability in an nonautonomous Lotka–Volterra system. *Nonlinear Analysis*, 40:37–49, 2000.

[BCL20] M.C. Bortolan, A.N. Carvalho, and J.A. Langa. *Attractors Under Autonomous and Non-autonomous Perturbations*, volume 246 of *Mathematical Surveys and Monographs*. AMS, 2020.

[BF20] F. Battelli and M. Fečkan. On the exponents of exponential dichotomies. *Mathematics*, 8:651, 2020.

[BP15] F. Battelli and K.J. Palmer. Criteria for exponential dichotomy for triangular systems. *Journal of Mathematical Analysis and Applications*, 428:525–543, 2015.

[CLR13] A.N. Carvalho, J.A. Langa, and J.C. Robinson. *Attractors for Infinite-dimensional Non-autonomous Dynamical Systems*, volume 182 of *Applied Mathematical Sciences*. Springer-Verlag, 2013.

[Cop65] W.A. Coppel. *Stability and Asymptotic Behaviour of Differential Equations*. Heath Mathematical Monographs. D.C.Heath, 1965.

[Cop78] W.A. Coppel. *Dichotomies in Stability Theory*, volume 629 of *Lecture Notes in Mathematics*. Springer-Verlag, 1978.

[CRLO23] A.N. Carvalho, L.R.N. Rocha, J.A. Langa, and R. Obaya. Structure of non-autonomous attractors for a class of diffusively coupled ODE. *Discrete and Continuous Dynamical Systems - Series B*, 28:426–448, 2023.

[CV02] V.V. Chepyzhov and M.I. Vishik. *Attractors for equations of mathematical physics*. AMS, 2002.

[Gop86a] K. Gopalsamy. Global asymptotic stability in a periodic Lotka–Volterra system. *The ANZIAM Journal*, 27:66–72, 1986.

[Gop86b] K. Gopalsamy. Global asymptotic stability in an almost periodic Lotka–Volterra system. *The ANZIAM Journal*, 27:346–360, 1986.

[Gue17] G.F. Guerrero. *Dinámica de redes mutualistas en ecosistemas complejos*. PhD thesis, Universidad de Sevilla, Sevilla, España, 6 2017.

[HS22] J. Hofbauer and S.J. Schreiber. Permanence via invasion graphs: incorporating community assembly into modern coexistence theory. *Journal of Mathematical Biology*, 85:Article number: 54, 2022.

[KR11] P. Kloeden and M. Rasmussen. *Nonautonomous Dynamical Systems*, volume 176 of *Mathematical Surveys and Monographs*. American Mathematical Society, 2011.

[Red96] R. Redheffer. Nonautonomous Lotka–Volterra systems. I. *J.Differential Equations*, 127:519–541, 1996.

[Tak96] Y. Takeuchi. *Global asymptotic dynamical properties of Lotka-Volterra systems*. World Scientific Publishing, 1996.